\documentclass[a4paper]{cas-sc}

\usepackage[authoryear,longnamesfirst]{natbib}
\usepackage{graphics}
\usepackage{graphicx}
\usepackage{mwe}
\usepackage{caption}
\usepackage{subcaption}
\usepackage{amsmath, amssymb, amsthm, amsfonts}
\usepackage{hyperref}
\usepackage{gensymb}
\usepackage[dvipsnames]{xcolor}
\usepackage{todonotes}
\usepackage{algorithm}
\usepackage{algpseudocode}
\usepackage{adjustbox}

\def\tsc#1{\csdef{#1}{\textsc{\lowercase{#1}}\xspace}}
\tsc{WGM}
\tsc{QE}
\tsc{EP}
\tsc{PMS}
\tsc{BEC}
\tsc{DE}

\begin{document}

\ExplSyntaxOn
\cs_gset:Npn \__first_footerline:
  { \group_begin: \small \sffamily \__short_authors: \group_end: }
\ExplSyntaxOff

\let\WriteBookmarks\relax
\def\floatpagepagefraction{1}
\def\textpagefraction{.001}
\shorttitle{}
\shortauthors{}

\title [mode = title]{Real-Time Calibration of Disaggregated  Traffic Demand}                      



\author[1]{Mozhgan Pourmoradnasseri}
[orcid=0000-0002-2092-816X]
\ead{mozhgan@ut.ee}
\cormark[1]

\author[1]{Kaveh Khoshkhah}
[ orcid=0000-0003-4777-804X]
\ead{kaveh.khoshkhah@ut.ee}


\author[1]{Amnir Hadachi}
[orcid=0000-0001-9257-3858   ]
\ead{hadachi@ut.ee}


\address[1]{ITS Lab, Institute of Computer Science, University of Tartu, Narva mnt 18, 51009 Tartu, Estonia}

\cortext[cor1]{Corresponding author}


\begin{abstract}  
Real-time traffic demand estimation is essential for intelligent transportation systems and traffic forecasting in urban areas. Hence, this paper presents a simulation-based optimization framework for city-scale real-time estimation and calibration of dynamic demand models by focusing on disaggregated microsimulation in congested networks. The calibration approach is based on sequential optimization demand estimation for short time frames and uses a stream of traffic count data from IoT sensors on selected roads. 
The proposed method builds upon the standard bi-level optimization formulation. The upper-level optimization problem is presented as a bounded variable quadratic programming in each time frame, making it computationally tractable. The DTA (dynamic traffic assignment) problem at the lower level is approached based on the stochastic route choice model. 
For every time frame, the probabilistic parameters of the route choice model are obtained through several rounds of feedback loop between the OD optimization problem (at the upper level) and parallel samplings and simulations for DTA (at the lower level). At the end of each time frame, the microsimulation state of the network is transferred to the next time frame to ensure the temporal dependency and continuity of the estimations during the time. In addition, the algorithm's convergence is accelerated by applying the fixed-point method to road-segment travel times.

The proposed sequential calibration model presents a drastic computational advantage over available methods by splitting the computations into short time frames and under the assumption that the demand in each time frame only depends on current and previous time frames. The model does not depend on reliable prior demand information. Moreover, with receiving the field measurements as a stream and the efficient time complexity of the algorithm in each time frame, the method successfully presents a solution for high-dimensional real-time and online applications. 
We validated the method with synthetic data and for a real-world case study in Tartu city, Estonia. The results show high accuracy under a tight computational budget.

\end{abstract}



\maketitle

\section{Introduction}

With the continued urbanization, commuters spend considerable time and energy in traffic, which imposes inevitable costs on the private and public sectors. To mitigate the losses, there is an increasing need for leveraging the capacity of existing infrastructure and improving the level of service through the advanced traffic management system. 

Real-time traffic monitoring systems are an essential component of intelligent transportation systems and decision support for analyzing and evaluating transportation facilities' current and future performance, observing and predicting the traffic condition, and measuring the short and long-term impacts of perturbation. However, due to the inherent complexity of urban dynamics, the accuracy, scalability, and efficiency of the solution approach remain challenging.
The smart city phenomenon and Internet of Things (IoT) data availability allowed the development of real-time traffic monitoring systems for intelligent transportation \citep{s22083030}. Many cities have taken this opportunity to collect and manage quality data for public and private use. The availability of a real-time data stream provides the possibility of developing detailed models for real-time estimation of mobility flow, even for light modes of transportation such as bicycles \citep{kumar2016inferring} and pedestrians \citep{kaveh2022}.

Data-driven approaches for dynamic demand modeling have been extensively studied during the last decades \citep{cascetta1993dynamic, zhou2007structural, osorio2019dynamic, zhang2021improving}. The objective of models is to generate outputs that describe the observed parameters of the network, such as link traffic counts, travel time, or average speed, through simulation tools. As powerful tools, traffic simulators attempt to replicate the traffic flow propagation by iteratively modifying the simulated parameters to surveillance data in different levels of aggregation. At the finest level, microscopic traffic simulators try to simulate the movements of all individual travelers in the network through several interactions between supply and demand models and adjusting calibration parameters based on real-world observation. 
Obtaining accurate simulation models for large-scale congested networks is a challenging task. Scalability, accuracy, computational efficiency, and consistency with real-world measurements are among the important aspects of traffic simulators. However, achieving reliable and robust simulation outcomes is highly challenging due to the high dimensionality of the problem and the necessity of reliable a priori information. In addition, for staying relevant to real-time applications, a stream of field measurements and a strict computational budget for calibration algorithms are crucial.

Dynamic traffic demand is often modeled as time-dependent origin-destination (OD) matrices by indicating the aggregated number of trips between selected origins and destinations in the transportation network for defined time frames.
In this approach, the whole time period of interest is split into smaller time frames, and a set of associated OD matrices is generated. 
Realistic dynamic demand models should take into account the individual behavior of travelers in route choice with respect to changes in the supply model, such as congestion and travel time. To address this challenge, the Dynamic OD matrix Estimation (DODE) problem implies solving the Dynamic Traffic Assignment (DTA) problem \citep{chiu2011dynamic, balakrishna2007offline}. 
The common approach in the literature for solving the DODE problem is the simulation-based bi-level optimization approach \citep{ashok2002estimation}. In the upper level, the discrepancy between the estimated and observed parameters in the OD granularity is minimized. The estimated OD matrix is verified in the lower-level problem by solving the DTA problem by heuristically assigning routes to OD trips, commonly through iterative stochastic simulations. Supply and demand models in bi-level optimization formulation interact with each other until the parameters of propagated traffic flow converge to the real field measurements. In each iteration, parameters have to be tuned using a calibration method to guarantee the accuracy and reliability of the outcome. 
The bi-level optimization problem implicitly models the nonlinear relation of supply and demand models. The DTA of the lower level problem constrains the optimization problem in the upper level, and only the optimum solution to the DTA problem is the feasible solution to the upper level problem \citep{andersen2022dynamic}. However, performing several DTAs for each optimization task demand considerable computational resources.

Researchers exploit various formulations and solutions both for the objective function of the optimization problem on the upper level and for tackling the DTA problem on the lower level. The optimization problem for OD calibration is commonly modeled as least square formulation by minimizing the distances between estimated and observed road counts, on the one hand, and the estimated and prior OD matrix, on the other hand, \citep{cascetta1993dynamic, bierlaire2004efficient, djukic2017modified}. This optimization problem is underdetermined for large-scale networks since the number of roads with traffic measurement is usually several orders of magnitude smaller than the number of all roads in the network. The dimension of the optimization problem is on the order of the number of OD pairs, which is considered high-dimensional for real-world case studies \citep{osorio2019high}. 

The optimization problem is mainly approached in the literature by general-purpose gradient-based solutions, such as Stochastic Perturbation Simultaneous Approximation (SPSA) \citep{ben2012dynamic, antoniou2015w}. The DTA problem is solved commonly through a simulator tool with a black-box strategy. In contrast, efficient solutions for high-dimensional transportation problems require designing model-specific algorithms.    
A metamodel for high-dimensional OD calibration is introduced in \citep{osorio2019high, osorio2019dynamic}. In this study, a robust and scalable model with high computational efficiency is achieved by embedding analytical and differentiable problem-specific structures into the algorithm. In this work, the simulation tool is used to derive analytical problems, not as a black-box. In the same line of research,  \citep{arora2021efficient} applies the analytical meta model to the static OD calibration problem by formulating the network model with a linear system of equations with a dimension of the same order as the number of roads in the network. With this technique, the authors gained computational efficiency and enhanced scalability. 

The simulation-based solutions to the DTA problem are shown to be more realistic in real-world applications \citep{ben2012dynamic, balakrishna2006off}, while other data-driven approaches, such as data of floating cars \citep{tsanakas2022d} or cellular network data \citep{sohn2008dynamic}, are exploited to model the route choice behavior of the travelers. Moreover, with technological advances in recent years, machine learning-based approaches are receiving increasing attention \citep{otkovic2013calibration, shafiei2018calibration, andersen2022dynamic}.

The dynamic demand estimation is addressed with two different strategies. The first formulation estimates the whole set of time-dependent OD matrices \textit{simultaneously} for the entire time period in one calibration round. The second formulation estimates and calibrates OD matrices \textit{sequentially} one at a time for each time frame. Most optimization-based studies have used the simultaneous approach, since the available sequential methods may lose accuracy due to giving up some information in the transition between time frames. While, on the other hand, in simultaneous estimation, all input data, including the traffic counts for the entire period, must be given to the model in advance, which makes it incapable of handling a stream of data for real-time or online demand calibration. Although the sequential technique has a high potential for being adapted to online calibration, it has received limited attention in the context of real-time demand estimation \citep{cascetta1993dynamic, ashok2002estimation, hu2017sequential, yang2018dynamic}.

From the perspective of spatiotemporal granularity, most presented works in the literature stay at OD-level and for a limited time period. Some works approached the real-time demand calibration problem \citep{bierlaire2004efficient, barcelo2015robust}, on small-scale networks and in an off-line context. Methods for achieving demand estimation at a higher spatial resolution are proposed in \citep{flotterod2009cadyts, flotterod2014disaggregate} for uncongested networks, and with a considerably high computational cost due to the high number of simulation runs required for convergence.   

The main objective of this work is to present a sequential optimization-based technique for the real-time disaggregated travel demand estimation and calibration problem. More precisely, we present an approach for dynamic trip-based route flow and OD matrix estimation. The calibration procedure takes as input the stream of traffic counts, aggregated in each time frame, and a normalized trip distribution between OD pairs. As opposed to common approaches, our method is not vulnerable to availability and accuracy of a prior OD demand (seed matrix). 
The bounded variable quadratic programming formulation of the calibration problem has a substantial computational advantage over the common Generalized Least Square (GLS) formulation by replacing the gradient-based solutions approach for continuous optimization with an efficient linear least-squares minimization.
Among the bottlenecks of sequential optimization approaches, are the running time of iteratively solving DTA through simulation, as a black-box and the assumption of independence of estimations in consecutive time frames. 
 To overcome the time complexity of DTA problem, an updated database of reasonable routes with their current travel times is maintained. The DTA problem in each time frame is iteratively solved based on probabilistic route assignment to OD trips in parallel microsimulations. Also, at the end of each time frame, the micro-level network state will be loaded to the next time frame using Simulation of Urban MObility (SUMO) \citep{lopez2018microscopic} tool for trip simulations. As a result, by transferring the traffic status from the current time frame to the following one, no information is lost, and a continuous microsimulation is generated for the entire period. In addition, Steffensen's fixed-point method is applied to link travel times in each calibration round to accelerate the convergence of iterative demand estimation. 

This paper is organized as follows. Section \ref{sec:problem} presents the high-level formulation of b-level optimization problem and our method to approach the solution. Section \ref{sec:arch} discusses the details of the algorithm steps for a fixed time frame.  Section \ref{sec:fixed-point} explains the Steffensen's fixed-point method applied to our algorithm. The method is validated by applying on synthetic networks in Section \ref{sec:val-syn} and for Tartu network in section \ref{sec:val-Tartu}. Finally, section \ref{sec:conclusion} concludes with summarizing the main outcomes of the work. 


\section{Problem formulation} \label{sec:problem}
Transportation networks are commonly modeled by a directed graph $\mathcal{G=(N,L)}$, when $\mathcal{N}$ is the set of nodes and $\mathcal{L}$ is the set of directed links. The set of nodes $\mathcal{N}$, usually correspond to junctions in the road network, and each link corresponds to a uni-directional road segment between two junctions. Additional attributes such as number of lanes, speed limit, length, and geolocation are associated to each link. A set of traffic counters is located on the subset of links $Q\subseteq \mathcal{L}$. Let $|Q|=q$ be the number of traffic counters.

The sets of origins and destinations are subsets of nodes $\mathcal{N}$.  A trip is identified by its origin and destination as the ordered pair $w=(i,j)\in W$, where $W$ is the set of all OD pairs.
We consider discrete-time frames of equal size $\Delta$ to approach the dynamic demand estimation. The $t$th time frame is denoted by $t$ when $t = 0,1, 2, \cdots, T-1$.

Optimization-based bi-level DODE problem is approached with two different formulation. 
In simultaneous framework, the problem is solved for all time frames, at once and the standard least square formulation of problem is as follows
\begin{equation} \label{eq:simult}
\begin{array}{ll@{}ll}
\text{min}  & \displaystyle \beta_1 \Vert \mathbf{c} - \mathbf{\Tilde c}\Vert_2 +  \beta_2 \Vert{\bf X} - {\bf \Tilde X}\Vert_2, & \\
\text{subject to} & \mathbf{c} = \text{DTA}({\bf X})
\end{array}
\end{equation}
when ${\bf \Tilde X}, {\bf X}\in \mathbb{R_+}^{|W|\times T}$ are prior and estimated OD matrices, and $\mathbf{\Tilde c}, \mathbf{c}  \in \mathbb{Z_+}^{q\times T}$ are real and estimated link measurements. The $\text{DTA}({\bf X})$ solves the dynamic traffic assignment problem, commonly using simulation tools and through a black-box approach.  The weights $\beta_1$ and $\beta_2$ are assigned based on confidence and the importance of each input data. The problem gets all field measurements $\mathbf{\Tilde c}$ and prior OD matrix ${\bf \Tilde X}$, in the beginning of the process and estimates a time-dependent OD matrix ${\bf X}$, at once. 

In the sequential formulation, the optimization problem is solved $T$ times, for $t=0, 1, 2, \cdots, T-1$

\begin{equation} \label{eq:seq}
\begin{array}{ll@{}ll}
\text{min}  & \displaystyle \beta_1 \Vert \mathbf{c} - \mathbf{\Tilde c_t}\Vert_2 +  \beta_2 \Vert{\bf X_t} - {\bf \Tilde X_t}\Vert_2, & \\
\text{subject to} & \mathbf{c} = \text{DTA}({\bf X})
\end{array}
\end{equation}
when  ${\bf \Tilde X_t}, {\bf X_t}\in \mathbb{R_+}^{|W|}$, and $\mathbf{c}, \mathbf{\Tilde c_t}  \in \mathbb{Z_+}^{q}$.

In this work, we exploit the sequential framework and build our method upon the formulation of equation \ref{eq:seq}. 
Our work differentiates from the previous works in several aspects. 
A significant difference between our method formulation with equation \ref{eq:seq} is that we do not constrain the optimization problem. The user equilibrium is approached through an iterative feedback loop of solving the optimization problem and performing DTA through stochastic route choice for verifying the optimization result.

The input requirement in our method is an initial probability distribution of trips between OD pairs instead of an initial OD matrix. More precisely, with OD probabilities, we indicate the chances of a random trip between different OD pairs. This input can be obtained more easily from available public data such as population distribution. Then, using the field measurement data and probabilistic parameters of the network, the probability distribution of OD trips is scaled up to an initial (seed) OD matrix.  

Our proposed method is based on the assumption that, for a fixed time frame and a given route, trips are generated every second based on a fixed Bernoulli probability distribution. The high-level idea of the proposed method is calculating these probability distributions for each route and time frame. Then, the best match to the input data is selected through sampling and microsimulation.  As a result, in addition to OD estimation, the method can estimate the disaggregated route-level demand in each time frame. 
The probabilistic sequential approach demonstrates a considerable computational advantage over available methods in the literature and is a step forward in real-time calibration. 

In contrast to other sequential approaches, the proposed solution approach considers the dependency of consecutive time frames by transferring the network traffic state from each time frame to the following, and therefore, no valuable information is discarded.   
Moreover, it allows receiving the field measurements as a data stream for $t=0, 1, 2, \cdots$, moving towards real-time and online applications. 

In section \ref{sec:arch}, the steps of our algorithm for a fixed time frame are explained in detail, and thus, the time index is omitted in the rest of this paper. 

\section{Methodology Architecture} \label{sec:arch}
\begin{figure}[h] 
\centering
\includegraphics[width=\textwidth]{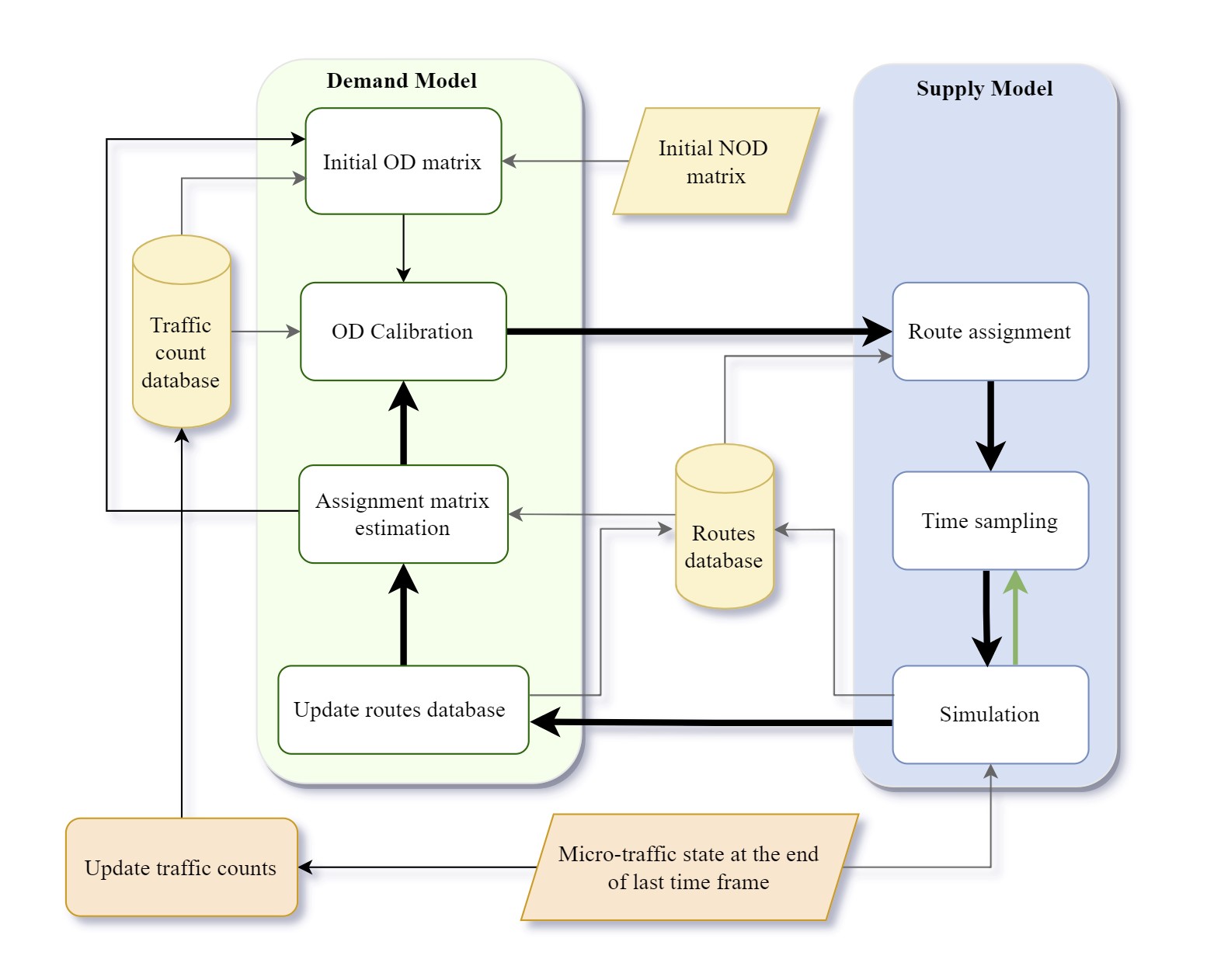}
\caption{High-level architecture of the calibration method.}
\label{fig:architecture}
\end{figure}

The high-level design of the method, for a fixed time frame $t$, is presented in figure \ref{fig:architecture}.  
In every time frame $t=0,1,2, \cdots$, the whole process will be performed in several iterations; therefore, the time index is omitted for simplicity. The exit condition is set in advance as an upper bound on the number of iterations of the whole process, or achieving a reasonably low error rate in the results of one iteration. 
In the following, we explain the elements of the method architecture.  

\subsection{Input data}
The algorithm depends on two input data; A rough estimation of trip distribution between OD pairs and link traffic counts. 
\subsubsection{Initial NOD matrix}
Our methodology relies on a prior demand estimation with minimal information. Let ${\bf N} = (\eta_{m})$ be a normalized matrix with $\sum_{m} \eta_{m} = 1$, such that each entry $\eta_{m}$ suggests an estimated ratio of trips between $m$th OD pair to all trips. Normalized trip distribution can be obtained from different data sources, such as travel surveys \citep{egu2020comparable}, census of population and housing \citep{arora2021efficient}, mobile data \citep{pourmoradnasseri2019od, hadachi2020unveiling}, public transport smart tickets \citep{mohamed2016clustering, yong2021mining}, or geo-tagged social media posts \citep{gao2014detecting}. 
We refer to this input as normalized OD (NOD) matrix. The standard decisive input for dynamic demand estimation is a reliable OD matrix, which is usually challenging to obtain. In contrast, the presented method relies only on trip weights between OD pairs and is not sensitive to the high accuracy of this prior estimation. 

\subsubsection{Traffic count database}
All dynamic models for traffic flow estimations require real-world measurement data for calibration. The more detailed and disaggregated this data, the more information can be extracted and better for parameter tuning. This work relies on disaggregated link traffic counts collected by stationary sensors. Moreover, for real-time applications, the presented method relies on the stream of traffic count data, possibly aggregated with a resolution of time frame length. 

\subsubsection{Routes database} \label{sec:route-database}
In the stochastic route choice model, OD trips are distributed to all routes between an OD pair. While enumerating all available routes imposes some theoretical and practical limitations. To reduce the complexity, we limit the method to a selection of reasonable routes between OD pairs and maintain these routes in a database. 

Let $R_{m}$ be the set of routes between  $m$th OD pair. Two approaches can be considered for creating the routes database. In the static approach, $R_{m}$ is extracted in advance from available resources and with their historical parameters for different time frames, such as Google Maps Distance Matrix API \citep{arora2021efficient}. In the dynamic approach, routes are updated gradually based on traffic status estimations in simulation rounds. In this work, we apply the dynamic approach. Initially, only the shortest paths for all OD pairs are added to the database. Then, in each iteration of the process and for every $m$, the best new route between OD pair $m$ with respect to travel time is added according to the latest simulation-based traffic status. The route database gradually grows and is maintained by storing the shortest paths between all OD pairs, with a predefined upper bound $\rho$ on the number of routes. Moreover, for each route $i\in R_{m}$, its travel time $\theta_{i}$ from the origin to destination and the travel times $\theta_{ik}$ from origin to sensors $k$ located on the route are computed and stored as route attributes. In order to be consistent, if a sensor $k$ is not located on a route $i$, we set $\theta_{ik}=\infty$. 

\subsection{Estimating assignment matrix}

The logit model is the common approach in the literature for route choice modeling \citep{osorio2019dynamic,arora2021efficient}. Let $R_{m}$ be the set of all routes (in the database) between OD pair $m$. The cost of route $i \in R_m$ is calculated by adding up the cost of all links belonging to the route, commonly indicated by  travel times $\theta_i$.  For the OD pair $m$, the trips are distributed among different routes based on logit stochastic route choice model. The probability $P_i$ of selecting route $i$ is proportional to the exponential of the cost of the route over the sum of all route costs, with some fixed parameter $\gamma$,
\begin{equation} \label{eq:logit}
    P_i=\frac{\exp ( \gamma \theta_i)}{\sum_{s\in R_m} \exp ( \gamma \theta_s)}. 
\end{equation}

In a selected time frame, we assume the trips departure time are drawn uniformly at random from the time interval. For a given route $i$, the probability $P_{ik}$ that a vehicle crosses a sensor $k$, depends on the length of time frame $\Delta$, and the travel time $\theta_{ik}$ from the origin to the location of sensor $k$,   
\begin{equation}\label{eq:prob}
 P_{ik} =     
    \begin{cases}
        \frac{\Delta-\theta_{ik}}{\Delta}, & \text{if}\ \Delta > \theta_{ik};\\
         0, & \text{otherwise}.
    \end{cases}
\end{equation}
Now, by merging equations \ref{eq:logit} and \ref{eq:prob}, he assignment matrix ${\bf A} = (\alpha_{mk}) $ is computed. Each matrix entry $\alpha_{mk}$ senotes the probability that a trip between $m$th OD pair crosses sensor $k$,

\begin{equation} \label{eq:sensor-crossig-pr}
    \alpha_{mk}=\sum_{i\in R_{m}} P_{ik} P_{i}.
\end{equation}

\subsection{Initial OD matrix}
The initial (seed) OD matrix ${\bf \Tilde{ X}}$, is computed by scaling NOD matrix $\bf N$, given as an input. The scaling factor is based on the sensor measurements and the latest assignment matrix probabilities of trips between different OD pairs, as follows. 
Let a vehicle plan $V_i^{\delta}$ be a trip on the route $i$, with the departure time of $\delta$.  We compute the scaling factor $\sigma$ to obtain an initial OD matrix $\bf\Tilde{X}$ from ${\bf N} = (\eta_{m})$ by

\begin{equation} \label{eq:scale}
    \sigma = \frac{\text{Sum of all sensor readings in the time frame}}{\text{Expected number of sensor hits by a random vehicle plan}}.    
\end{equation}
Let $s$ be the expected number of sensor hits by a random vehicle plan $V_i^{\delta}$ drawn from the set of all possible routes $i$ and with a random departure time $\delta$. Then, 
\begin{equation}
 s= \sum_{m} \eta_{m} \sum_{k} \alpha_{mk}.   
\end{equation}
The equation \ref{eq:scale} can be rewritten as
\begin{equation}
    \sigma = \frac{\sum_k {\Tilde c_k}}{s}, 
\end{equation}
with ${\Tilde c_k}$ being the traffic counts of sensor $k$, aggregated in the selected time frame. 
Finally, the initial OD-matrix $\bf\Tilde{X}$ is calculated by the scalar product of the scaling factor and the OD-matrix of trip distribution; 
\begin{equation}
  {\bf\Tilde{X}} = \sigma \cdot {\bf N}.  
\end{equation}

\subsection{OD calibration}
The calibration step aims to estimate the OD demand with minimum distance to the prior OD matrix and real-world measurements after route assignment. Let  ${\bf X}\in \mathbb{R_+}^{|W|}$ be the estimated OD matrix, with positive entries\footnote{For convenience, vector representation for matrices is used. }. Each entry $x_{m} \in {\bf X}$ shows the (expected) number of trips between $m$th OD pair, in the current time frame. In the context of the bi-level optimization problem, the upper level optimization problem is formulated as the following least square problem 
\begin{equation} \label{eq:first-opt}
\begin{array}{ll@{}ll}
\text{min}  & \displaystyle \beta_1 \Vert f({\bf X})\Vert_2 +  \beta_2 \Vert g({\bf X})\Vert_2, & \\
\end{array}
\end{equation}
 when $\Vert f({\bf X})\Vert_2 = \Vert{\bf A}{\bf X} - \mathbf{\Tilde c}\Vert_2 $ estimates the discrepancy between the expected traffic counts and the real field measurements on the links with sensors, and $\Vert g({\bf X})\Vert_2 = \Vert{\bf X} - {\bf \Tilde X}\Vert_2 $ estimates the discrepancy between the computed OD matrix and the prior OD matrix ${\bf \Tilde X}$ in the same time frame. The vector of observed link traffic counts is denoted by $\mathbf{\Tilde c}  \in \mathbb{Z_+}^{q}$. The weights $\beta_1$ and $\beta_2$ are assigned based on the importance of each input data.
 Equation \ref{eq:first-opt} can be formulated as follows
\begin{equation} \label{eq:second-opt}
\begin{array}{ll@{}ll}
\displaystyle  \min_{l\leq {\bf X} \leq u}  & \displaystyle \sqrt {\beta_1 \Vert f({\bf X})\Vert^2_2 + \beta_2 \Vert g({\bf X})\Vert^2_2 }& \\
\end{array}
\end{equation}
and by defining  ${\bf A'}= \begin{vmatrix}{\bf A} \\ \lambda{\bf I}\end{vmatrix} $ and ${\bf \Tilde b} = \begin{vmatrix}\mathbf{\Tilde c} \\ \lambda{\bf \Tilde X} \end{vmatrix}$, and bounding the demand, we obtain 
\begin{equation} \label{eq:third-opt}
\begin{array}{ll@{}ll}
\displaystyle  \min_{l\leq {\bf X} \leq u}  & \displaystyle  \Vert {\bf A'} {\bf X} -  {\bf \Tilde b} \Vert_2, & \\
\end{array}
\end{equation}
where $\lambda$ is a constant coefficient, reflecting the degree of confidence on prior NOD ${\bf N}$.

The optimization problem \ref{eq:third-opt}, introduces a bounded variable non-negative least square (NNLSQ) minimization, with large and sparse coefficient matrix. Several fast methods for solving NNLSQ exist, such as STIR \citep{branch1999subspace}, and BVLS \citep{stark1995bounded}. 
In this work, \texttt{scipy.optimize.lsq\_linear} solver \citep{NNLSQ} is effectively exploited.

\subsection{Route assignment} \label{sec:route-assignment}
The solution of the optimization problem \ref{eq:third-opt}, gives the real-valued matrix ${\bf X^*} \in \mathbb{R_+}^{|W|}$ indicating the expected number of trips between OD pairs. In this step,  the trips are distributed between different routes proportional to probabilities of equation \ref{eq:logit}.  
As a result, an aggregated route-level OD matrix is obtained. In other words, for every OD pair $m$ and route $i \in R_m$,
\begin{equation} \label{eq:route assignment}
E_i=x^*_m \times P_i   
\end{equation}
gives the real-valued expected number of trips in route $i$. 

\subsection{Time sampling and simulation} \label{sec:simulation}
This part of the algorithm provides a complementary step for the DTA problem. In order to disaggregate the route-level trips in equation \ref{eq:route assignment}, every traveler (agent) chooses a departure time $\delta$, based on the plan $V_i^{\delta}$, and all trips are simulated using the SUMO simulator. More precisely, in every time unit (one second), a trip is generated according to the Bernoulli distribution, with probability $E_i/\Delta$, on route $i$. With this approach, we avoid the rounding error for converting the number of trips from real to integer. 
The result of time sampling gives a set of vehicle plans $V_i^{\delta}$ that SUMO will simulate.

The main obstacle in performing simulation-based DODE is the computational cost associated with solving several DTA problems for evaluating the solution of the optimization problem \ref{eq:first-opt}. 
We tackle this problem first by probabilistic approach in route assignment (section \ref{sec:route-assignment}) and then by performing parallel time sampling and simulation. Several time samplings and simulations can be performed in parallel depending on the available resources and applications. Then, the best simulation result with respect to the discrepancy between the observed traffic count and the estimated traffic count can be selected. Moreover, the new best routes for all OD pairs, together with the travel times $\theta_{i}$ and $\theta_{ik}$, are added to the route database as it is discussed in section \ref{sec:route-database}. 
At the end of the current time frame, several vehicles may not have arrived at their destinations. We transfer the set of these uncompleted trips from the current time frame to the next one. SUMO simulator provides the possibility of saving and loading the network state. We use this option to save the network state at the end best simulation of the current time frame and reload the state at the beginning of the next time frame. It allows the vehicles that did not finish their trip during the current time frame to continue the trip at the beginning of the next time frame. As a result, although the problem is solved sequentially, time frames are not independent, and the output quality has all the advantages of simultaneous simulation.  

The main output of this step is considered to be the average travel time of each link in the network, according to the best time sampling and simulation. We denote this output by ${\bf \tau} \in \mathbb{R}_+^{|\mathcal{L}|}$ when $\mathcal{L}$ is the set of links in the network. 

\subsection{Updating route database}
In this step, in addition to updating the best routes between OD pairs, the probabilities associated with routes are updated. Link travel times from the simulation are converted to route probabilities based on the logit model in equation \ref{eq:logit}. 
The shortest routes between OD pairs, are selected using NetworX \footnote{\href{https://networkx.org/}{https://networkx.org/}} python package, by assigning the link travel times as weights to edges. 

\subsection{Updating traffic counts} 
This work assumes that the maximum trip duration between any OD pair is not longer than the time frame $\Delta$. Therefore, if a vehicle departs at time frame $t$, it has to arrive at its destination in any of the time frames $t$ or $t+1$. Thus, all incomplete trips loaded from $t$ to $t+1$ will terminate their trips in $t+1$.  Nevertheless, this set of trips may influence the traffic counts in $t+1$ by crossing the sensor locations in the rest of their trips, while being accounted for in the demand estimation of time frame $t$. To solve this issue, the expected number of sensor hits caused by incomplete trips in $t$ will be removed from the traffic count data of $t+1$.

\section{Enhancement by fixed-point method approach} \label{sec:fixed-point}
The bi-level optimization formulation of demand estimation in a congested network can be seen as a fixed-point problem. The upper level problem estimates an OD matrix that, when assigned to the network at the lower level, regenerates the network parameters used in estimating the OD matrix at the upper level. In practice, within iterative interactions between the supply model and the demand model, the stochastic user equilibrium is achieved \citep{bell1997transportation}. Different fixed-point methods are applied on link demand vector \citep{meng2004transportation, cascetta2001fixed},for the set of links with traffic counters. In this work we apply Steffensen's method on link travel time vector, for all network links. 

\begin{figure}[h] 
\centering
\includegraphics[width=.8\textwidth]{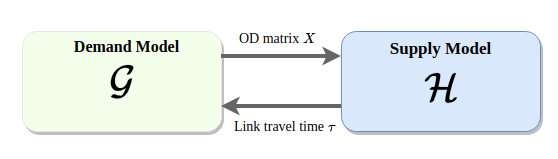}
\caption{Fixed-point schematic}
\label{fig:fixed-point}
\end{figure}

 We can present our demand model by function
\begin{equation}
  \mathcal{G}: \mathbb{R}^{|\mathcal{L}|} \rightarrow \mathbb{R}^{|W|} \quad \text{with} \quad
  \mathcal{G}(\tau) = X,
\end{equation}
and the supply model by function
\begin{equation}
  \mathcal{H}: \mathbb{R}^{|W|} \rightarrow \mathbb{R}^{|\mathcal{L}|} \quad \text{with} \quad
  \mathcal{H}(X) = \tau.
\end{equation}
In the other words, as it is presented in figure \ref{fig:fixed-point}, the demand model $\mathcal{G}$ estimates the number of trips between OD pairs, $X$, based on the estimated link travel times, $\tau$. Then, given the OD matrix $X$, the supply model calculates the travel times of each link. 

In each round of iteration, we measure the distance between the current link travel times and the previous travel time vector $\tau$. A smaller distance indicates the consistency of the solution. In other words, by defining $\mathcal{F=H\circ G}$, the problem is translated to a fixed-point problem by searching for the link travel time vector $\tau$, such that $\mathcal{F}(\tau) = \tau$. 

Common technique for estimating the solution of fixed-point problems is through iterative approaches that $\tau_{n+1}=\mathcal{F}(\tau_{n})$. Two main problems with iterative techniques are slow convergence and the possibility of divergence. In the stochastic route choice model, an increase in the travel time results in a decrease in the probabilities of assignment matrix $A$ (equation \ref{eq:sensor-crossig-pr}). Consequently, to satisfy sensor readings, the estimated number of trips increases; hence, in the supply model, travel times increase further. Repeating this process may lead to divergence from the fixed-point. To prevent this issue, we bound the travel time for each link. If $\nu$ is the vector of free flow travel time of different links, we set $\nu \leq \tau \leq d \cdot \nu$, when $d$ is a suitable constant. 
Then, the revised version of $\mathcal{F}$ is defined as 
\begin{equation}
    \bar{\mathcal{F}}(\tau) = \max\{ \min\{\mathcal{F}(\tau), d\cdot\nu\}, \nu \}.
\end{equation}

To accelerate the method's convergence, we use Aitken's delta-squared process \citep{suli2003introduction} in Steffensen's root-finding method \citep{johnson1968steffensen}, presented in Algorithm \ref{alg:steffensen}. The main advantage of Steffensen's method over other fixed-point methods with quadratic convergence is replacing derivative calculation with the first-order divided difference of the function.

\begin{algorithm}
\caption{Fixed-point algorithm based on Steffensen's method}\label{alg:steffensen}
    \hspace*{\algorithmicindent} \textbf{Input} \text{$\tau_0$, IterationNumber}\\
    \hspace*{\algorithmicindent} \textbf{Output} $\tau^*$
\begin{algorithmic}
\State $N \gets 0$
\State $\text{minError} \gets \infty$
\While{$N < \text{IterationNumber}$}
    \State $N \gets N + 1$
    \State $\tau_1 \gets \bar{\mathcal{F}}(\tau_0)$
    \State $\tau_2 \gets \bar{\mathcal{F}}(\tau_1)$
    \State $\tau \gets \tau_0 - {(\tau_1 - \tau_0)^2}/{(\tau_2 - 2\tau_1 +\tau_0})$
    \State $\tau \gets \max\{ \min\{\tau, d\cdot\nu\}, \nu \}$ \Comment{To force $\nu \leq \tau \leq d \cdot \nu$. }
    \State $error \gets \text{RelataiveError}(\tau_0, \tau)$
    \If{$error < \text{minError}$}
        \State $\text{minError} \gets error$
        \State $\tau^* \gets \tau$
    \EndIf
    \State $\tau_0 \gets \tau$
\EndWhile
\end{algorithmic}
\end{algorithm}


\section{Validation setup}
In order to evaluate the result of the proposed method, we conducted tests on two synthetic networks, and for the case-study of Tartu city. 
To quantify the approximations' accuracy, the common root mean square error (RMSE), normalized root-mean-square error (NRMSE), and the relative error ($\epsilon$) are used:

\begin{equation}
   \text{RMSE} ({\bf Y}, {\bf Y}^*) =\frac{1}{\sqrt{N}}\Vert{\bf Y} - {\bf Y}^*\Vert_2 ,
\end{equation}

\begin{equation}
   \text{NRMSE} ({\bf Y}, {\bf Y}^*) =\frac{\text{RMSE} ({\bf Y}, {\bf Y}^*)}{{\bf \bar{y}}},
\end{equation}
and
\begin{equation}
   \epsilon({\bf Y}, {\bf Y}^*) = \frac{\Vert{\bf Y} - {\bf Y}^*\Vert_2 }{\Vert{\bf Y}\Vert_2 } \times 100, 
\end{equation}
where ${\bf Y}$ and ${\bf Y}^*$ are the vector of true and estimated variables and $N$ is the number of variables to be estimated, and ${\bf \bar{y}}$ is the mean value of the entries of vector ${\bf Y}$.

All experiments in this work is performed on a machine with an Intel Core i7-9800X sixteen-core processor and 64 GB of RAM, running on Ubuntu. 

\section{Validation with synthetic network} \label{sec:val-syn}In the first test, a smaller network  and synthetic traffic data is produced for time period of $4$ hours. Then the synthetic input data is given to the algorithm  to replicate the observed flow on network links. In the end, the performance of the algorithm is measured based on the estimated link traffic counts in consecutive time frames. 

In the second experiment, a bigger network is generated and for one fixed time frame the robustness of algorithm and its convergence is analysed through different iterations of algorithm and various error measurements. 

\subsection{Small Network}
In this section we first explain the details of synthetic network and traffic data generation, and then we discuss the algorithm performance on the given data. 
\subsubsection{Synthetic data generation}
\begin{figure}[h] 
\centering
\includegraphics[width=.3\textwidth]{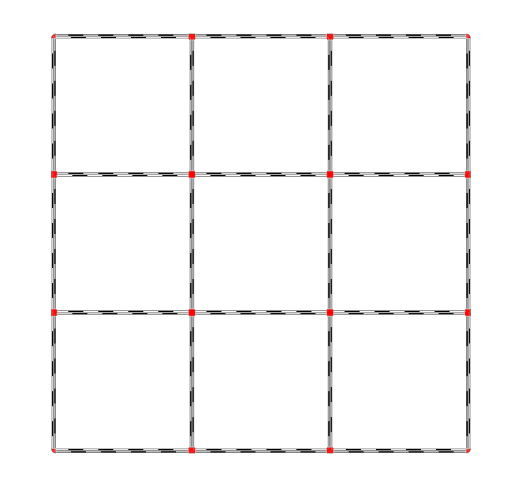}
\caption{The synthetic network used in the experiment. All links are equipped with a traffic sensor and all junctions have a traffic light}
\label{fig:network}
\end{figure}

A synthetic grid network is designed as figure \ref{fig:network}, with $16$ nodes and $48$ directed links between them. Each presented link is $400$ meters long and contains two lanes. Every junction is equipped with a traffic light with protocols embedded from SUMO, and every link with a traffic counter.

All junctions are considered as points of interests, and therefore, in the presented grid, $240$ OD pairs are available. 
We generated a set of random trips between every OD pairs, with a random number of trips chosen uniformly between $200$ and $350$ in the time interval of $4$ hours. Then, by running duarouter \footnote{\href{https://sumo.dlr.de/docs/duarouter.html}{https://sumo.dlr.de/docs/duarouter.html}} process in SUMO for $100$ times, routes are assigned to trips. The best set of trips, with respect to travel time is selected and simulated for collecting synthetic traffic count data. After this process, a total number of $62,500$ trips are generated. 

We split the time interval, into $4$ one-hour time frames. For each time frame, link traffic counts are aggregated and stored as sensor data for the input of the algorithm . For the second input, NOD, we normalize the OD matrix of generated trips in the whole $4$ hours. 

\subsubsection{Results}
The synthetic sensor data and NOD are given as inputs to the algorithm, and the estimated link counts after calibration are compared against the input (ground truth). To observe the effect of parameter tuning in accuracy of the output, we run the experiment for three different values of parameter $\lambda=1$, $\lambda= 0.01$ and $\lambda= 100$, in the optimization problem \ref{eq:third-opt}. $\lambda$ reflects the weight of the prior OD matrix in the objective function of the calibration model. The more confident we are about the initial OD matrix, we can increase the value of $\lambda$. In our method, we only use trip distribution as external information, and the initial OD is calculated by  the scaling up parameter $\sigma$. Therefore, as observed in the results, for smaller values of $\lambda$, a better match to ground truth traffic count is achieved, and for larger values of $\lambda$, the distance between estimated and initial OD matrices decreases. 

The presented results are obtained after $100$ rounds of calibration, with the execution of $30$ time sampling and simulation in parallel, in every round. Each iteration of the algorithm is performed on an average of $45$ seconds. As a result, the total running time of the method for each one-hour time frame adds up to approximately $75$ minutes. 

For $\lambda=1$, figure \ref{fig:traffic counts, lambra=1} shows the comparison of ground truth and estimated link traffic counts for different time frames. Table \ref{table:error lambda=1} presents different error measures of the proposed method at each time for OD matrix and link traffic counts. In the designed network, since all links are equipped with sensors, the error calculations take into account all network links. 
\begin{figure}[pos=htb!]
     \centering
     \begin{subfigure}[b]{0.245\textwidth}
         \centering
         \includegraphics[width=\textwidth]{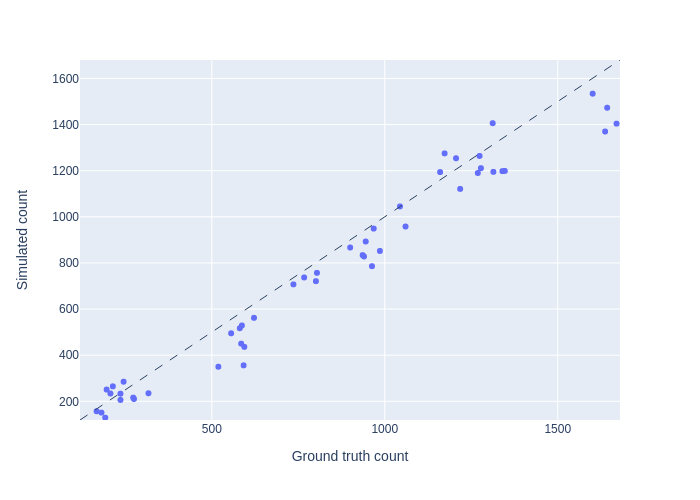}
         \caption{Time frame $t=1$}
         \label{fig:y equals x}
     \end{subfigure}
     \hfill
     \begin{subfigure}[b]{0.245\textwidth}
         \centering
         \includegraphics[width=\textwidth]{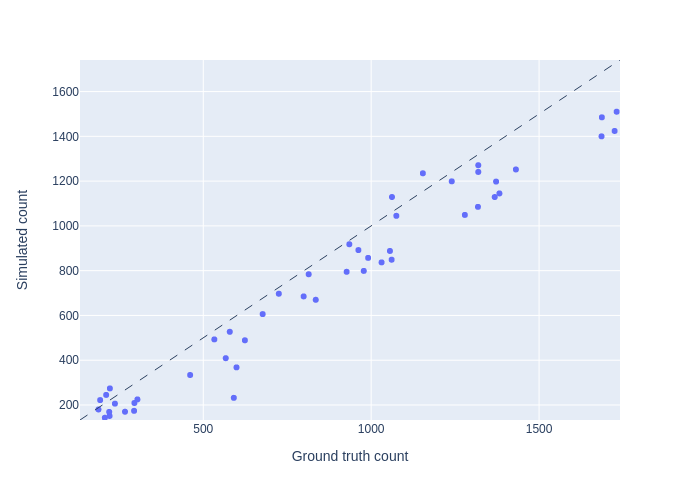}
         \caption{Time frame $t=2$}
         \label{fig:three sin x}
     \end{subfigure}
     \begin{subfigure}[b]{0.245\textwidth}
         \centering
         \includegraphics[width=\textwidth]{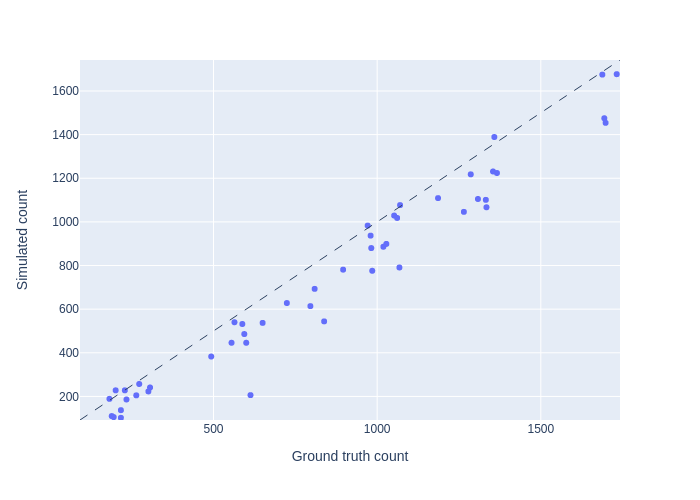}
         \caption{Time frame $t=3$}
         \label{fig:five over x}
     \end{subfigure}
          \hfill
     \begin{subfigure}[b]{0.245\textwidth}
         \centering
         \includegraphics[width=\textwidth]{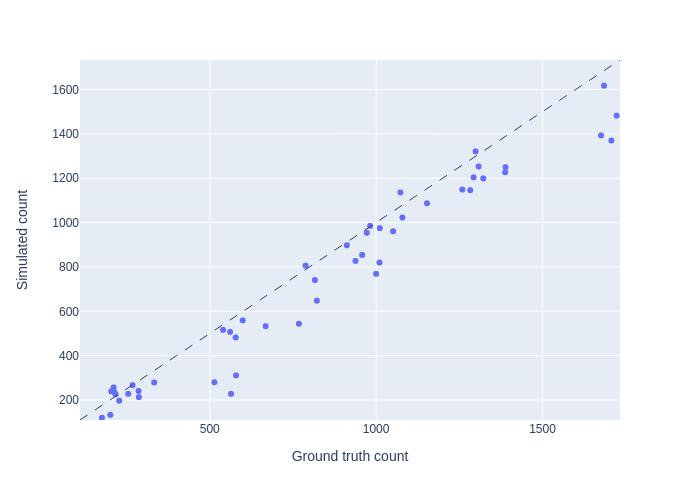}
         \caption{Time frame $t=4$}
         \label{fig:three sin x}
     \end{subfigure}
        \caption{Ground truth counts and simulated counts in different time frames for $\lambda=1$. }
        \label{fig:traffic counts, lambra=1}
\end{figure}
\begin{table}[pos=htb!]
\centering
\caption{Trip counts and error estimations in different time frames for $\lambda=1$}
\begin{adjustbox}{max width=\textwidth}
\begin{tabular}{c|ll|lll|lll}
Time frame & Real count & Estimated count & OD $\epsilon$ \% & OD RMSR & OD NRMSE \% & Sensor $\epsilon$ \% & Sensor RMSR & Sensor NRMSE \% \\ \hline
0          & 16074      & 15586           & 82.42            & 56.39   & 84.2        & 11.28                & 106.26      & 12.85           \\
1          & 15550      & 14843           & 111.91           & 74.59   & 115.12      & 15.58                & 151.23      & 17.81           \\
2          & 15495      & 13172           & 86.04            & 57.17   & 88.56       & 14.95                & 143.53      & 17.07           \\
3          & 15379      & 13831           & 108.2            & 71.5    & 111.58      & 14.15                & 135.03      & 16.14           \\ 
\end{tabular}\label{table:error lambda=1}
\end{adjustbox}
\end{table}

Figure \ref{fig:traffic counts, lambda=0.01} and table \ref{table:error lambda=0.01} show the results and errors 
for $\lambda=0.01$, and different time frames. 
\begin{figure}[pos=h]
     \centering
     \begin{subfigure}[b]{0.245\textwidth}
         \centering
         \includegraphics[width=\textwidth]{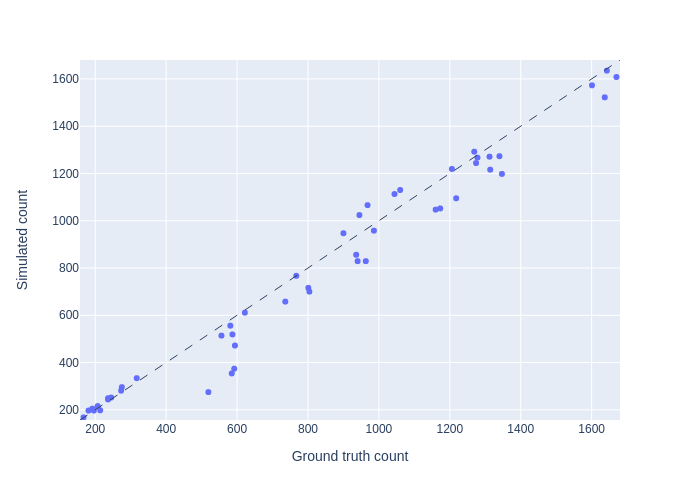}
         \caption{Time frame $t=1$}
         \label{fig:y equals x}
     \end{subfigure}
     \hfill
     \begin{subfigure}[b]{0.245\textwidth}
         \centering
         \includegraphics[width=\textwidth]{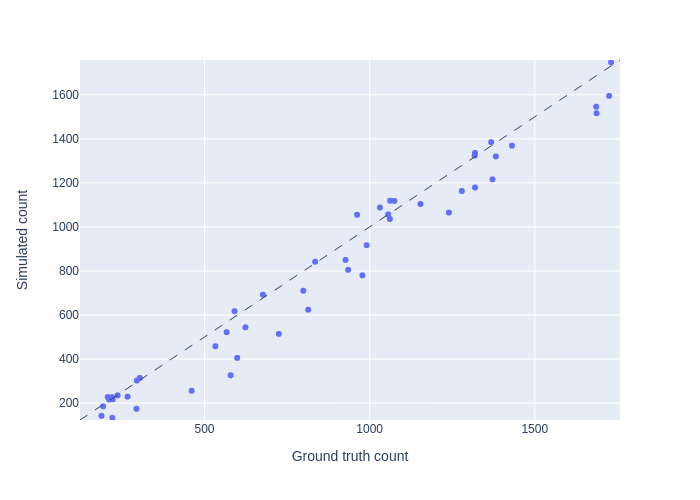}
         \caption{Time frame $t=2$}
         \label{fig:three sin x}
     \end{subfigure}
     \begin{subfigure}[b]{0.245\textwidth}
         \centering
         \includegraphics[width=\textwidth]{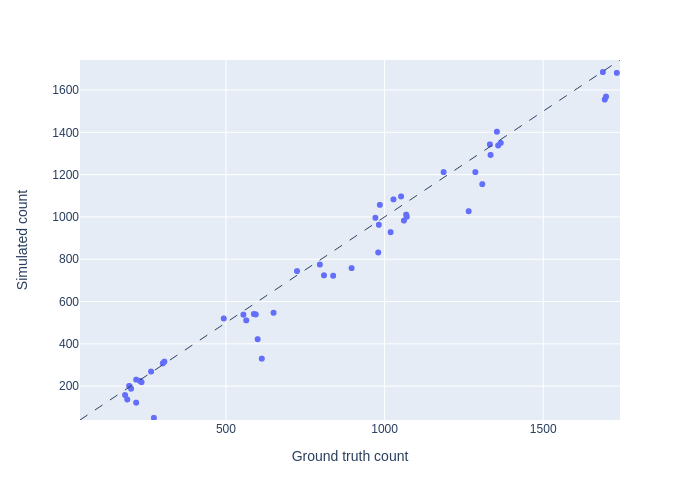}
         \caption{Time frame $t=3$}
         \label{fig:five over x}
     \end{subfigure}
          \hfill
     \begin{subfigure}[b]{0.245\textwidth}
         \centering
         \includegraphics[width=\textwidth]{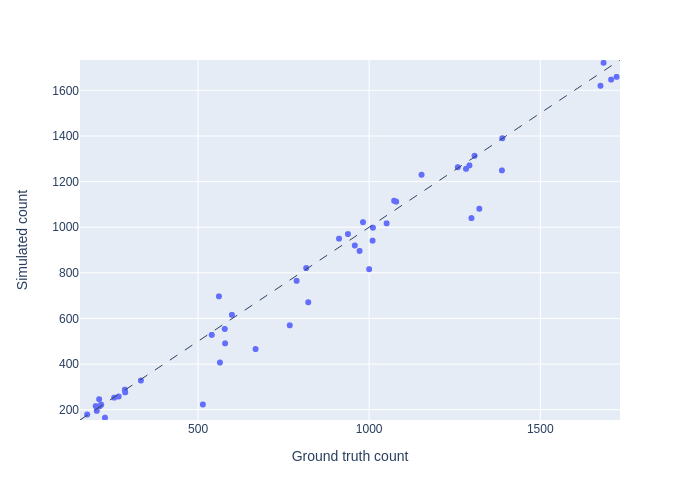}
         \caption{Time frame $t=4$}
         \label{fig:three sin x}
     \end{subfigure}
        \caption{Ground truth counts and simulated counts in different time frames, for $\lambda=0.01$. }
        \label{fig:traffic counts, lambda=0.01}
\end{figure}
\begin{table}[pos=h]
\centering
\caption{Trip counts and error estimations in different time frames for $\lambda=0.01$}
\begin{adjustbox}{max width=\textwidth}
\begin{tabular}{c|ll|lll|lll}
Time frame & Real count & Estimated count & OD $\epsilon$ \% & OD RMSR & OD NRMSE \% & Sensor $\epsilon$ \% & Sensor RMSR & Sensor NRMSE \% \\ \hline
0          & 16074      & 17097           & 126              & 86.22   & 128.73      & 9.36                 & 88.23       & 10.67           \\
1          & 15550      & 17867           & 199.77           & 133.15  & 205.51      & 10.76                & 104.44      & 12.3            \\
2          & 15495      & 15076           & 120.11           & 79.81   & 123.61      & 9.74                 & 93.54       & 11.12           \\
3          & 15379      & 16051           & 172.67           & 114.1   & 178.06      & 10.3                 & 98.28       & 11.75           \\ 
\end{tabular} \label{table:error lambda=0.01}
\end{adjustbox}
\end{table}

The results of the method in different time frames and for $\lambda=100$ are presented in figure \ref{fig:traffic counts, lambda=100}, and table \ref{table:error lambda=100}. 
\pagebreak

\begin{figure}[pos=htb!]
     \centering
     \begin{subfigure}[b]{0.245\textwidth}
         \centering
         \includegraphics[width=\textwidth]{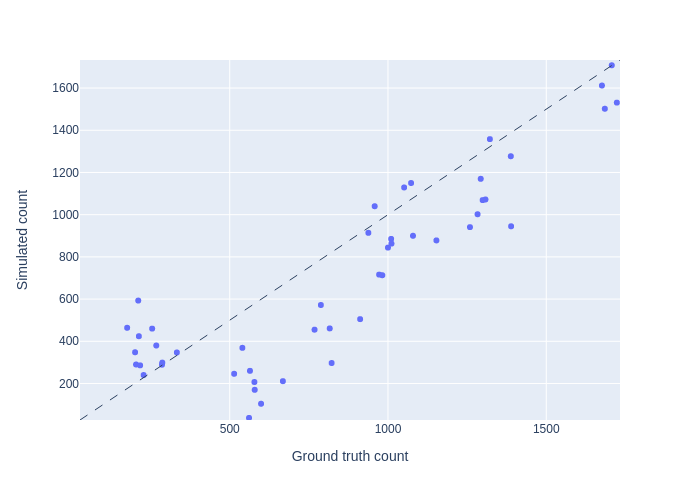}
         \caption{Time frame $t=1$}
         \label{fig:y equals x}
     \end{subfigure}
     \hfill
     \begin{subfigure}[b]{0.245\textwidth}
         \centering
         \includegraphics[width=\textwidth]{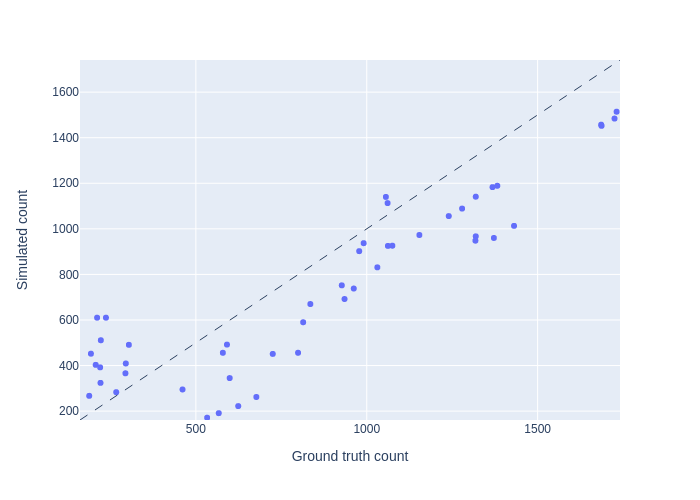}
         \caption{Time frame $t=2$}
         \label{fig:three sin x}
     \end{subfigure}
     \begin{subfigure}[b]{0.245\textwidth}
         \centering
         \includegraphics[width=\textwidth]{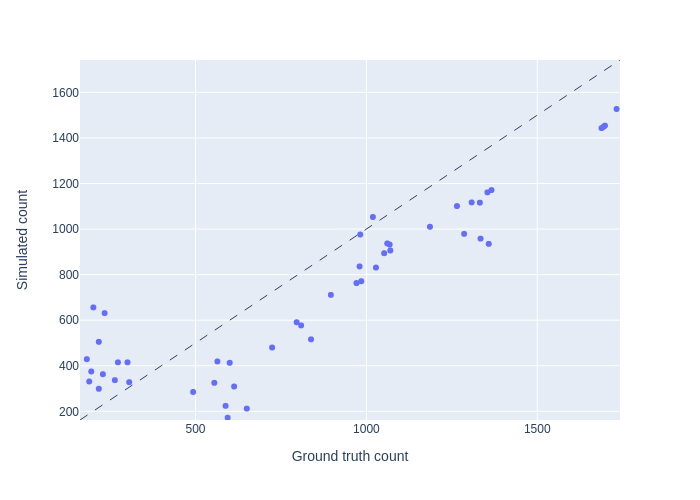}
         \caption{Time frame $t=3$}
         \label{fig:five over x}
     \end{subfigure}
          \hfill
     \begin{subfigure}[b]{0.245\textwidth}
         \centering
         \includegraphics[width=\textwidth]{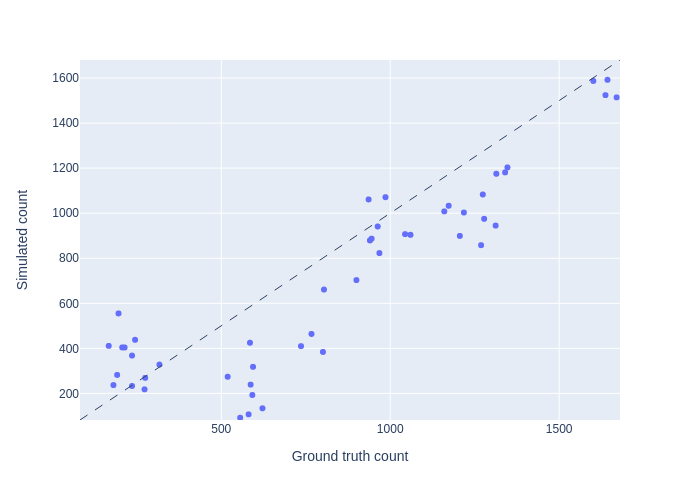}
         \caption{Time frame $t=4$}
         \label{fig:three sin x}
     \end{subfigure}
        \caption{Ground truth counts and simulated counts in different time frames, for $\lambda=100$. }
        \label{fig:traffic counts, lambda=100}
\end{figure}

\begin{table}[pos=htb!]
\centering
\caption{Trip counts and error estimations in different time frames for $\lambda=100$}
\begin{adjustbox}{max width=\textwidth}
\begin{tabular}{c|ll|lll|lll}
Time frame & Real count & Estimated count & OD $\epsilon$ \% & OD RMSR & OD NRMSE \% & Sensor $\epsilon$ \% & Sensor RMSR & Sensor NRMSE \% \\ \hline
0          & 16074      & 12893           & 19.81            & 13.56   & 20.24       & 25.17                & 237.11      & 28.69           \\
1          & 15550      & 13386           & 13.93            & 9.29    & 14.34       & 25.07                & 243.4       & 28.66           \\
2          & 15495      & 13256           & 14.48            & 9.62    & 14.91       & 24.94                & 239.49      & 28.49           \\
3          & 15379      & 12485           & 18.82            & 12.43   & 19.4        & 27.13                & 258.87      & 30.95           \\ 
\end{tabular} \label{table:error lambda=100}
\end{adjustbox}
\end{table}

\subsection{Large Network}
In this section the experiment is run on a larger network and a more in-depth analysis of robustness of the algorithm is discussed, in a fixed time frame. 

\subsubsection{Synthetic data generation}
\begin{figure}[h] 
\centering
\includegraphics[width=.3\textwidth]{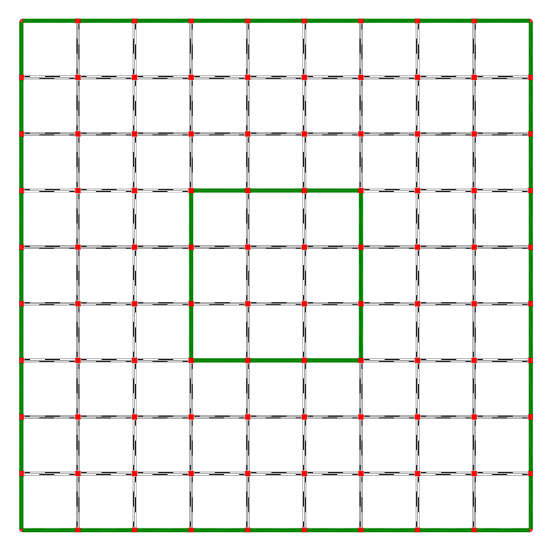}
\caption{The synthetic network used in the experiment. All links are equipped with a traffic sensor and all junctions have traffic lights. Nodes located on green edges are points of interest as the potential origins and destinations of trips.}
\label{fig:10-network}
\end{figure}

The $10\times 10$ synthetic  grid network is generated similar to the previous case. The grid consists of $100$ nodes and $360$ directed links between them. Each presented link is $200$ meters long and contains two lanes. Every junction is equipped with a traffic light with protocols embedded from SUMO, and every link with a traffic counter. The nodes of the outer ring and an inner ring, which are distinguished by green in figure \ref{fig:10-network}, are selected as points of interest. All nodes are considered as potential origins/destinations. In the end, random trips are generated with a rate between $0$ and $15$ for all $2256$ OD pairs in the duration of one hour. Then, similar to the previous case, route assignment is performed after $100$ runs of duarouter, and selecting the best set of simulated trips with respect to travel time. As a result, a total number of $16829$ trips are generated and simulated and link traffic counts are collected.

\subsubsection{Results}
After $150$ iterations of the algorithm, a relative error of $\epsilon=8\%$ in the estimation of link traffic counts is achieved.
In order to analyze the robustness of the presented method, we investigate different error measurements over different iterations of the algorithm. 

First, we remind and define a few parameters for one iteration round of calibration and simulation. 
${\bf \Tilde{c}}$ is the vector of ground truth traffic counts on network links, given as the input. 
${\bf A}$ is the assignment matrix (as defined in equation \ref{eq:sensor-crossig-pr}), and ${\bf X^*}$ is the estimated OD matrix after calibration (solution of optimization problem \ref{eq:third-opt}).
Let ${\bf c^*}$ be the vector of estimated link traffic counts after route assignment and simulation at the end of one iteration (obtained at the end of section \ref{sec:simulation}). Let ${\bf \tau^*}$ be the vector of experienced link average travel times after the best time sampling and simulation at the end of the current iteration round, as described in section \ref{sec:simulation}. Denote the average speed of the whole network in the current iteration by $\bar s$.

In each iteration round of calibration and simulation, the following errors are propagated in different steps of the algorithm: 
\begin{enumerate}
    \item {\bf OD calibration error $\epsilon ({\bf  {\Tilde c}, AX^*})$} is the relative error between ground truth link traffic count data ${\bf \Tilde{c}}$ and the expected link traffic counts ${\bf AX^*}$, after OD calibration step. 
    \item {\bf Calibration-to-simulation error $\epsilon ({\bf  AX^*, c^*})$} is the relative error between the expected link traffic count data, ${\bf AX^*}$, obtained from OD calibration and the vector of estimated link traffic counts, ${\bf c^*}$ after sampling and simulation step. 
    \item {\bf Iteration error $\epsilon ({\bf  {\Tilde c}, c^*)}$ } is the relative error between ground truth link traffic count data ${\bf \Tilde{c}}$ and approximated link traffic count ${\bf c^*}$ at the end of current iteration. \item {\bf Fixed-point error $\epsilon ({\bf \tau_t}, {\bf \tau_{t+1}})$} is the relative error between average travel time vector ${\bf \tau}$ in the current and previous iterations. 
\end{enumerate}

Table \ref{table:covariance} compares the joint variability of the errors observed in different iterations of the process. 
Iteration error measures the distance between field traffic counts on each link and the estimated link traffic counts, and it is commonly used in the literature as a metric for quantifying the performance of demand estimation. 

The final objective of our algorithm is achieving a network state estimation with low error in one of the several iteration rounds. Then, the output of that estimation will be selected as the procedure's output. One iteration round compounds two main steps of OD calibration and route assignment and simulation, each producing some errors that affect the iteration error. The high correlation of iteration error is observable with OD calibration error and calibration-to-simulation error. In the calibration phase, the OD is estimated based on the available route choice probabilities computed in previous iterations. Therefore, inaccurate estimation of the mentioned parameters results in error propagation in the entire process. 

Time and route sampling in the route assignment and simulation step also may lead to error generation. Moreover, the inaccurate route travel times in the route database will affect the simulation output negatively. Link travel time estimations also influence fixed-point errors. Therefore, as expected, we observe a high correlation between fixed-point and calibration-to-simulation error, while the correlation is moderate in the case of OD calibration.

\begin{table}[pos=h]
\centering
\caption{Covariance between error values and the average speed in different iteration rounds.}
\begin{adjustbox}{max width=\textwidth}
\begin{tabular}{l| c c c c c} 

 {\bf Covariance} & OD calibration error & Calibration-to-simulation error & Iteration error & Fixed-point error & Average speed \\ [1ex] 
 \hline \\
OD calibration error &   1 &  0.733 & 0.964 & 0.632 & -0.585 \\ [1ex]
Calibration-to-simulation error &  0.733  & 1 & 0.884 &  0.820 & -0.879 \\ [1ex]
Iteration error &  0.964 & 0.884 & 1 & 0.741 &  -0.737\\[1ex]
Fixed-point error &  0.632 &  0.820 &0.741 & 1 & -0.827\\[1ex]
Average speed & -0.585  &  -0.879 & -0.737 & -0.827 & 1\\ [1ex] 

\end{tabular}
\end{adjustbox}
\label{table:covariance}
\end{table}

\begin{figure}[pos=h]
     \centering
     \begin{subfigure}[b]{0.49\textwidth}
         \centering
         \includegraphics[width=\textwidth]{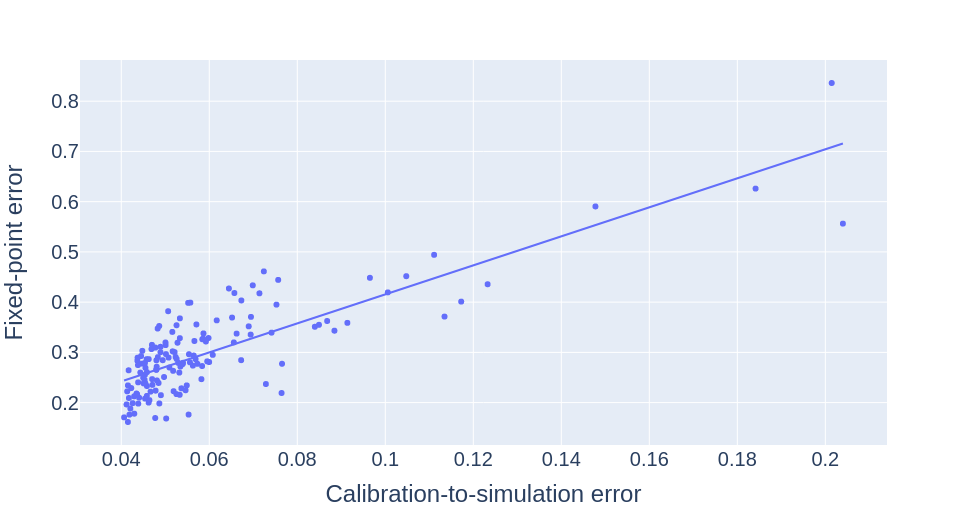}
         \caption{Relation between fixed-point error and calibration-to-simulation error.}
         \label{fig:y equals x}
     \end{subfigure}
     \hfill
     \begin{subfigure}[b]{0.49\textwidth}
         \centering
         \includegraphics[width=\textwidth]{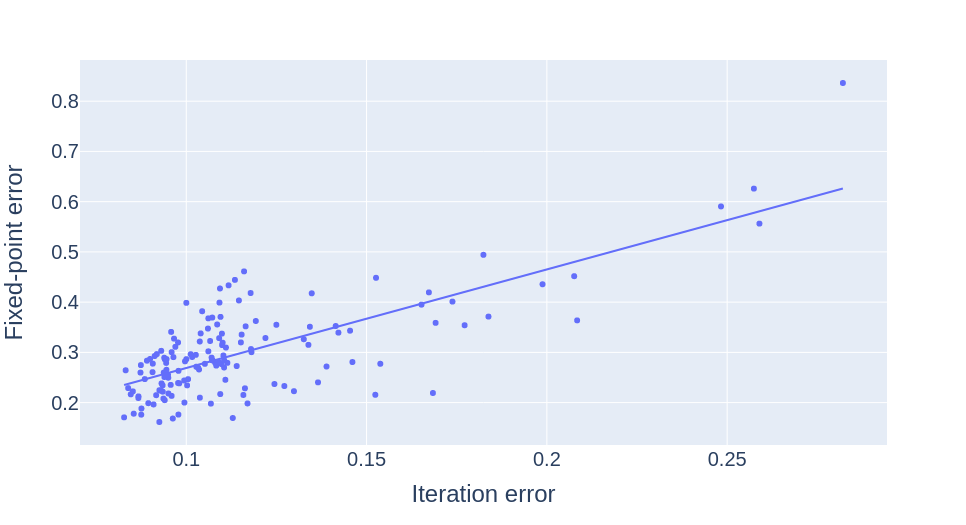}
         \caption{Relation between fixed-point error and iteration error.\break}
         \label{fig:three sin x}
     \end{subfigure}
        \caption{Comparison of fixed-point error with calibration-to-simulation error and iteration error in different iterations of the algorithm.  }
        \label{fig:fixed-point-error}
\end{figure}

\begin{figure}[pos=h]
     \centering
     \begin{subfigure}[b]{0.49\textwidth}
         \centering
         \includegraphics[width=\textwidth]{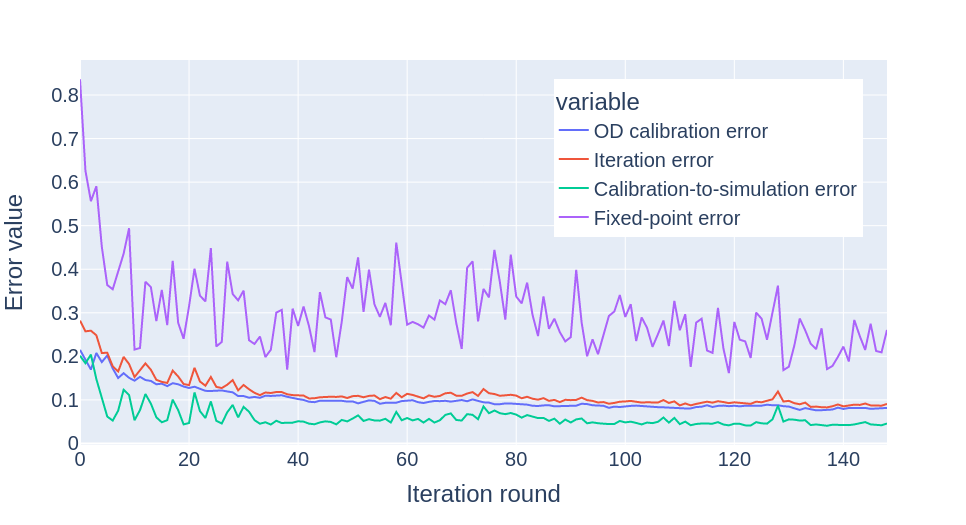}
         \caption{Error values in each iteration of the algorithm. \break}
         \label{fig:all-errors}
     \end{subfigure}
     \hfill
     \begin{subfigure}[b]{0.49\textwidth}
         \centering
         \includegraphics[width=\textwidth]{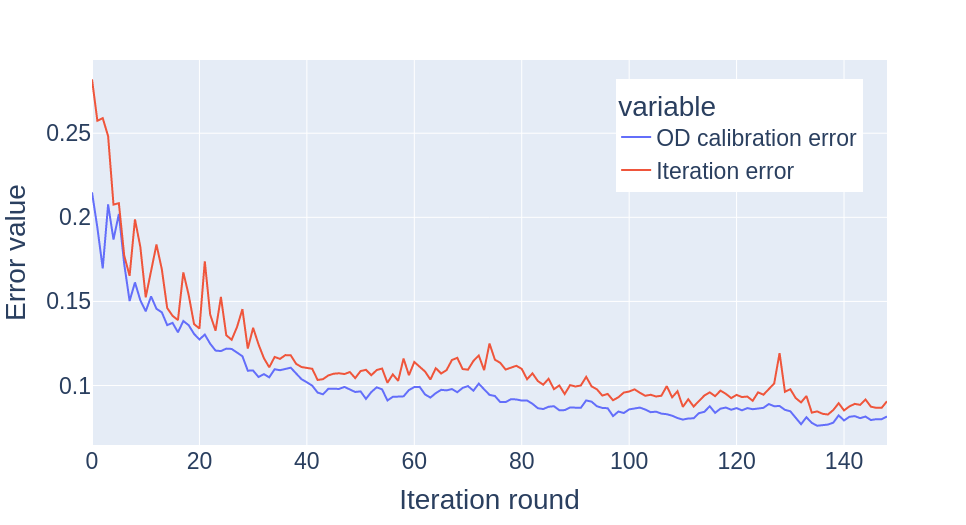}
         \caption{OD calibration and iteration error values in different iterations.}
        \label{fig:calibration-iteraton-error}
     \end{subfigure}

        \caption{Error trends in different iterations of the algorithm.  }
\end{figure}

As discussed in section \ref{sec:fixed-point} and algorithm \ref{alg:steffensen}, in order to accelerate the convergence of the algorithm, we aim to minimize fixed-point error as the main metric of the performance of the method. Figure \ref{fig:fixed-point-error} shows the relation of fixed-point method with calibration-to-simulation error and iteration error, in different iterations of the algorithm. 
Figure \ref{fig:all-errors} presents the variability of the discussed errors in $150$ iteration rounds, and figure \ref{fig:calibration-iteraton-error} focuses on OD calibration and iteration errors with a higher resolution. With more iterations of the algorithm, more routes are added to route database and we observe a decreasing trend in the defined errors.

\begin{figure}[pos=h] 
\centering
\includegraphics[width=.5\textwidth]{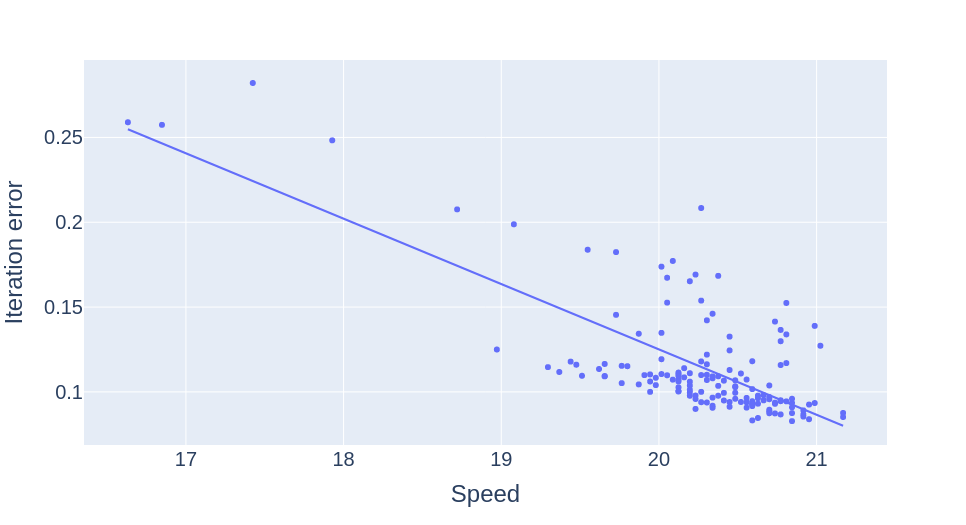}
\caption{The relation between average speed and iteration error. }
\label{fig:speed}
\end{figure}

A high negative correlation can be observed between average speed, $\bar s$, with fixed-point error, and calibration-to-simulation error. Figure \ref{fig:speed} shows the relation between the iteration error and average speed in each round of the algorithm. Minimizing the network's average speed can also be an alternative approach for achieving the equilibrium. 

\section{Case Study: Tartu} \label{sec:val-Tartu}
\begin{figure}[pos=h] 
\centering
\includegraphics[width=.8\textwidth]{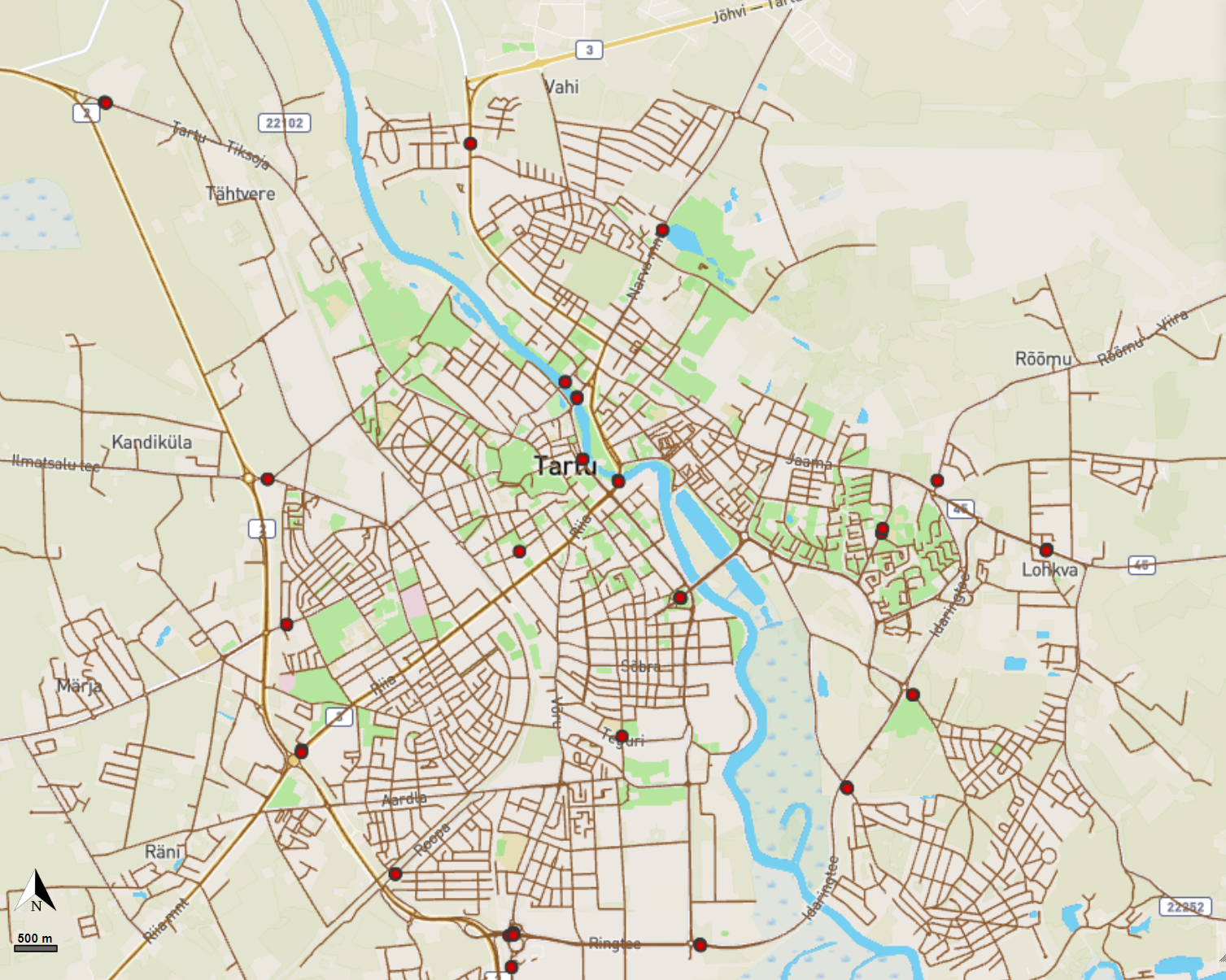}
\caption{Road network of the city of Tartu with the locations of IoT traffic count sensors.}
\label{fig:tartu}
\end{figure}

The performance of the proposed method in a city-scale network is evaluated in the case study of the city of Tartu, the second-largest city in Estonia. The city's approximate land area and population are $39~km^2$ and $91000$.
The road network model of the study area consists of $2780$ nodes and $6700$ links. There are $33$ IoT traffic counters located in the city. 
The AVC sensors placed on the city's border give us a big picture of the vehicles entering and leaving the city of Tartu. In addition, Thinnect sensors provide us with a glimpse into the mobility pattern of inner-city areas. The data is transferred to the Cumulocity platform, managed by the city, in real-time. Figure \ref{fig:tartu} presents the road network of the city and the location of IoT traffic counters.

A set of origins/destinations are selected proportional to the number of network nodes in each city district. As a result, $28820$ OD pairs are considered in this case study.
We run the method for a full day and sequentially estimate $24$ OD matrices, each corresponding to a one-hour time frame. In addition, the probabilistic approach with sampling and simulation and using the SUMO simulator generates a continuous microsimulation of all individual trips for the whole $24$h. 

\subsection{Experiment setup}
The method requires two input datasets; traffic counts and the prior NOD trip distribution. As mentioned before, the Cumulocity platform provides us with access to the stream of traffic count data. The prior NOD matrix can be obtained from reliable sources such as census data. In this study, we generated the data according to district populations and network topology. 
The city consists of $17$ administrative districts. In addition, to consider the commuters across the city borders, an auxiliary district is defined to reflect intercity trips. The trip distribution between inside city (in) and outside city (out) is based on following weights: in/in = 0.65, in/out = 0.15, out/in = 0.15, and out/out = 0.05. First, an $18\times 18$ NOD is generated on the district level, proportional to the population of each district and according to the defined weights. Then, trips between districts are distributed uniformly between OD pairs inside a district.

Finally, the method is fed by the hourly traffic count data and prior NOD matrix. The method design aims for efficiency under a tight computational budget. 
Therefore, the exit condition is selected as an upper bound of 10 calibration rounds or achieving a relative error of less than $10\%$ for link count estimates (iteration error $\epsilon(\bf \Tilde{c} , c^*)$). In each calibration round, we perform time sampling and simulation steps 15 times in parallel, and the best result with respect to the ground truth data is transmitted to the next calibration round. Although the solution approach is based on sequential simulation optimization, SUMO simulation provides the possibility of loading the micro-level network state at the end of each time frame to the following one. Therefore, in the end, a continuous microsimulation output for the whole period is also generated.  

\subsection{Results}
\begin{figure}[pos=h]
     \centering
     \begin{subfigure}[b]{0.49\textwidth}
         \centering
         \includegraphics[width=\textwidth]{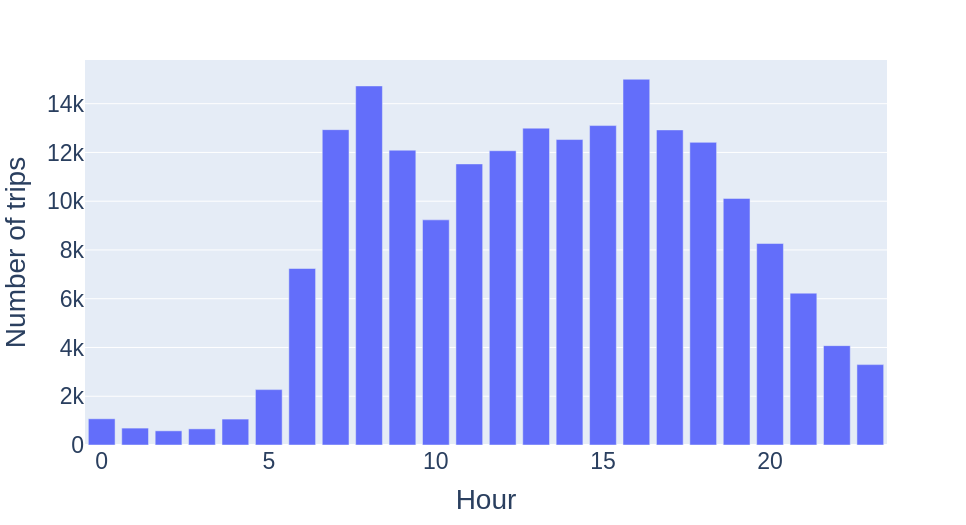}
         \caption{Distribution of number of trips in different hours.}
         \label{fig:trip number}
     \end{subfigure}
     \hfill
     \begin{subfigure}[b]{0.49\textwidth}
         \centering
         \includegraphics[width=\textwidth]{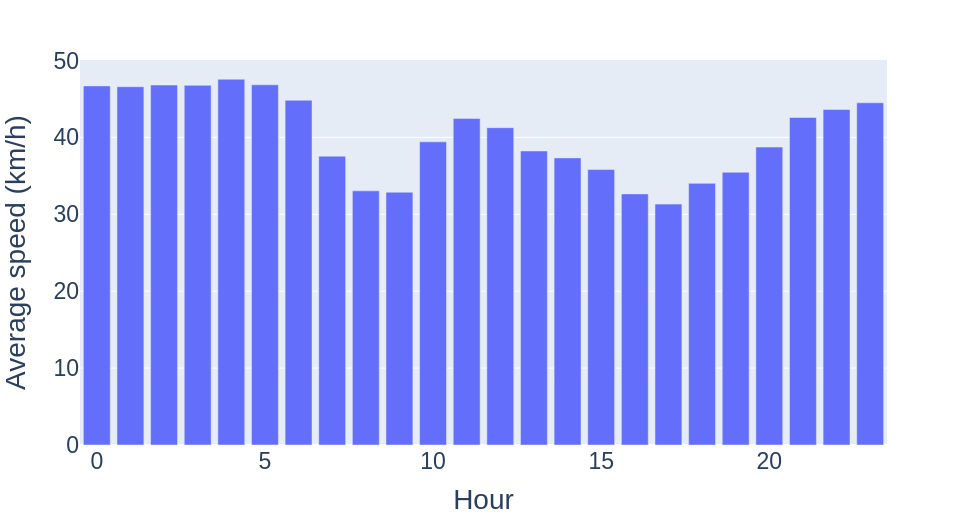}
         \caption{Average speed in different hours.}
         \label{fig:trip speed}
     \end{subfigure}

        \caption{Estimated number of trips and estimated average speed in different hours of the selected day. }
        \label{fig:trips}
\end{figure}

In this section, the results of the case study are presented. Figure \ref{fig:trip number} shows the hourly distribution of the estimated number of trips. Figure \ref{fig:trip speed} shows the hourly distribution of the network's average speed of the simulated trips, after the method execution. Speed limits in Tartu city are generally 50 km/h or 30 km/h in built-up areas. The correlation between the traffic density and average speed is observed in figure \ref{fig:trips}. 

\begin{figure}[h!]
     \centering
     \begin{subfigure}[b]{0.23\textwidth}
         \centering
         \includegraphics[width=\textwidth]{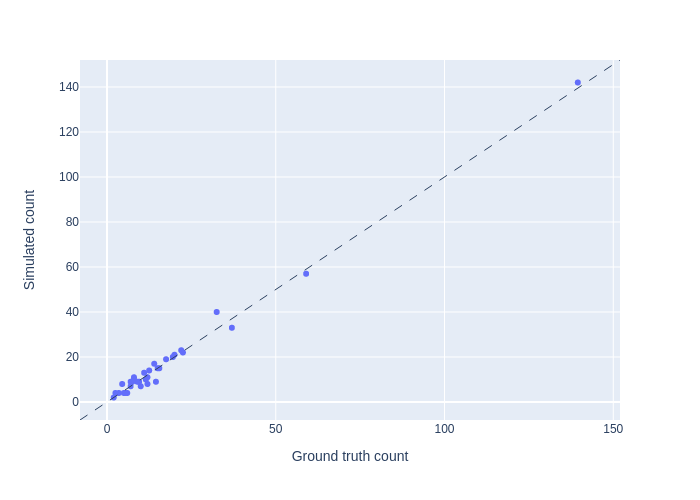}
         \caption{00:00 to 01:00}
         \label{fig:y equals x}
     \end{subfigure}
     \hfill
    \begin{subfigure}[b]{0.23\textwidth}
         \centering
         \includegraphics[width=\textwidth]{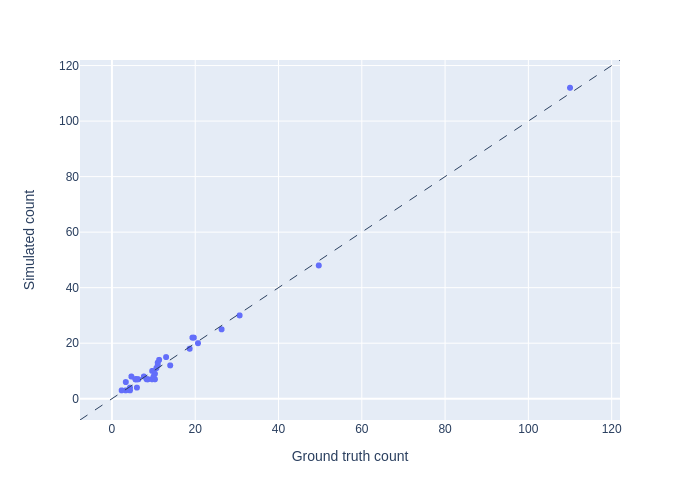}
         \caption{01:00 to 02:00}
         \label{fig:y equals x}
    \end{subfigure}
     \hfill
    \begin{subfigure}[b]{0.23\textwidth}
         \centering
         \includegraphics[width=\textwidth]{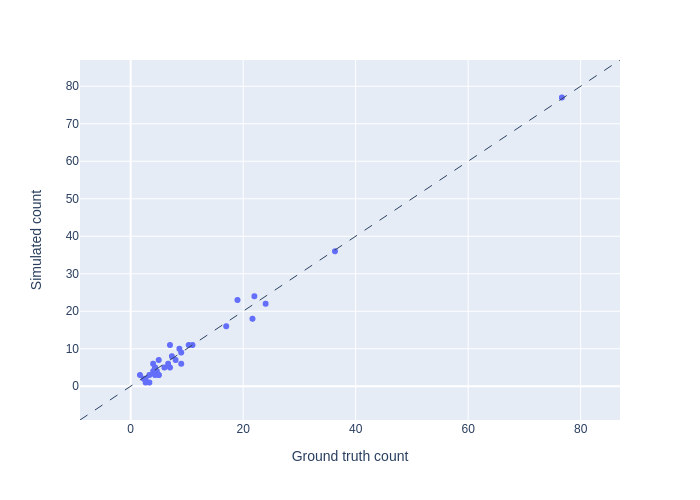}
         \caption{02:00 to 03:00}
         \label{fig:y equals x}
    \end{subfigure}
     \hfill
     \begin{subfigure}[b]{0.23\textwidth}
         \centering
         \includegraphics[width=\textwidth]{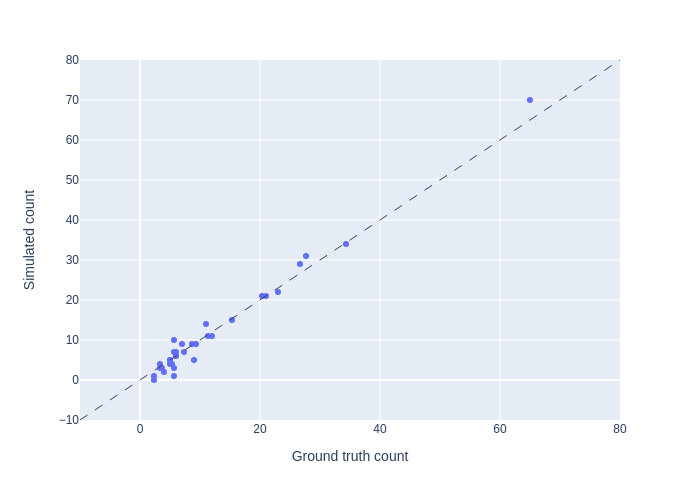}
         \caption{03:00 to 04:00}
         \label{fig:three sin x}
     \end{subfigure}
     \\
     \begin{subfigure}[b]{0.23\textwidth}
         \centering
         \includegraphics[width=\textwidth]{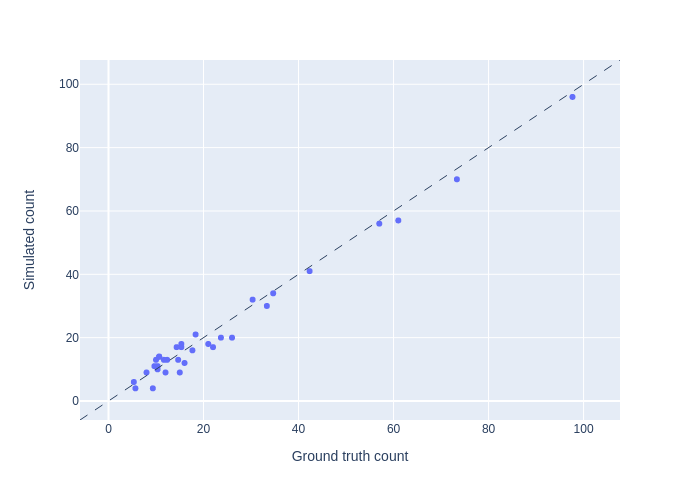}
         \caption{04:00 to 05:00}
         \label{fig:y equals x}
     \end{subfigure}
     \hfill
    \begin{subfigure}[b]{0.23\textwidth}
         \centering
         \includegraphics[width=\textwidth]{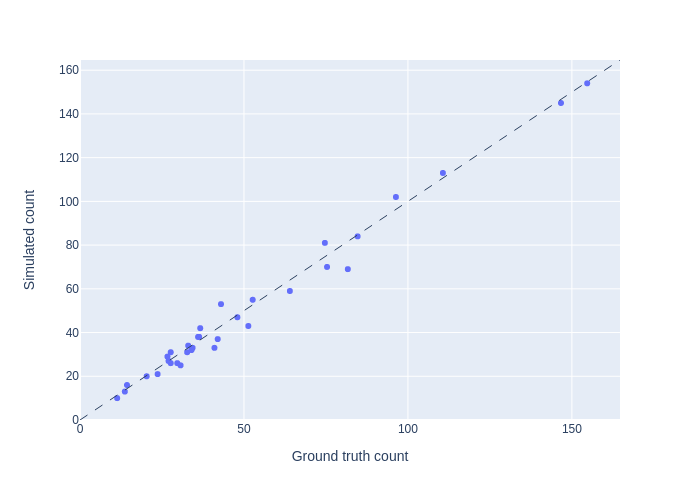}
         \caption{05:00 to 06:00}
         \label{fig:y equals x}
    \end{subfigure}
     \hfill
    \begin{subfigure}[b]{0.23\textwidth}
         \centering
         \includegraphics[width=\textwidth]{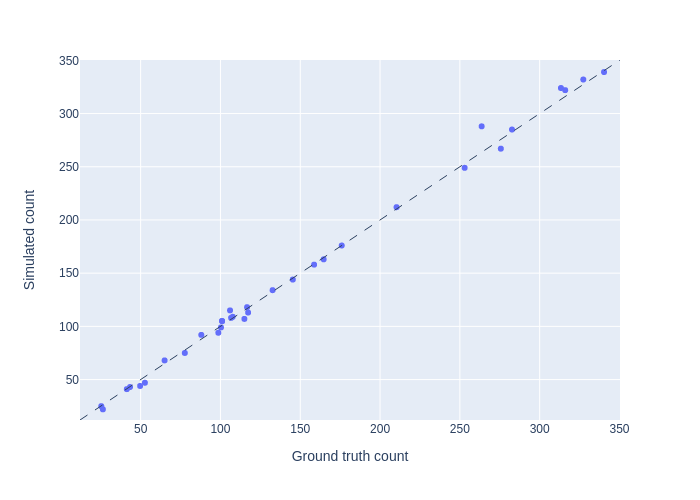}
         \caption{06:00 to 07:00}
         \label{fig:y equals x}
    \end{subfigure}
     \hfill
     \begin{subfigure}[b]{0.23\textwidth}
         \centering
         \includegraphics[width=\textwidth]{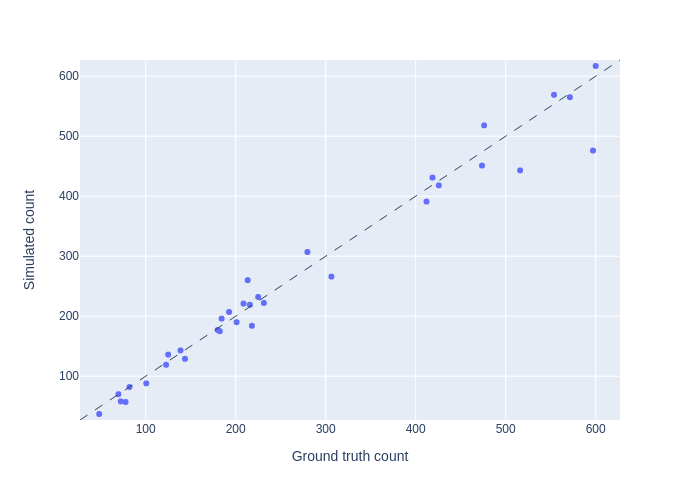}
         \caption{07:00 to 08:00}
         \label{fig:three sin x}
     \end{subfigure}
     \\
     \begin{subfigure}[b]{0.23\textwidth}
         \centering
         \includegraphics[width=\textwidth]{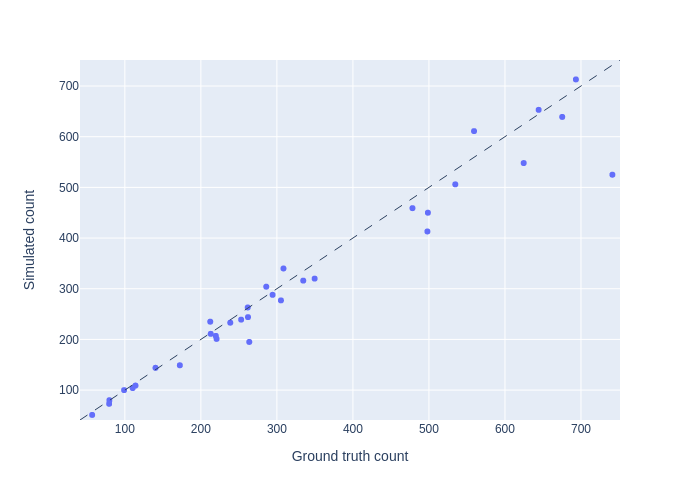}
         \caption{08:00 to 09:00}
         \label{fig:y equals x}
     \end{subfigure}
     \hfill
    \begin{subfigure}[b]{0.23\textwidth}
         \centering
         \includegraphics[width=\textwidth]{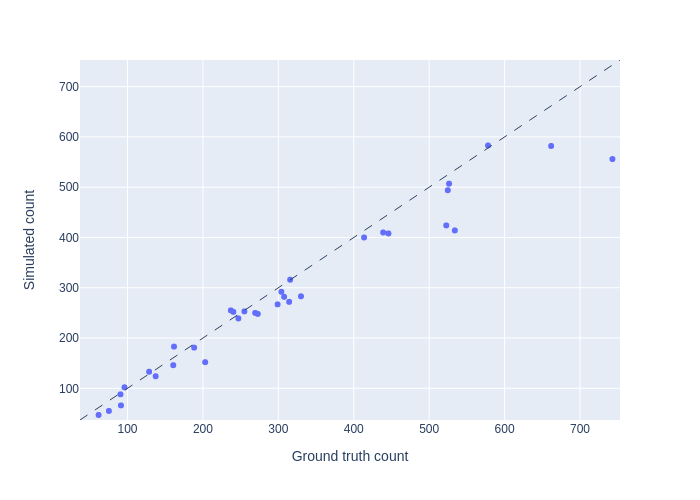}
         \caption{09:00 to 10:00}
         \label{fig:y equals x}
    \end{subfigure}
     \hfill
    \begin{subfigure}[b]{0.23\textwidth}
         \centering
         \includegraphics[width=\textwidth]{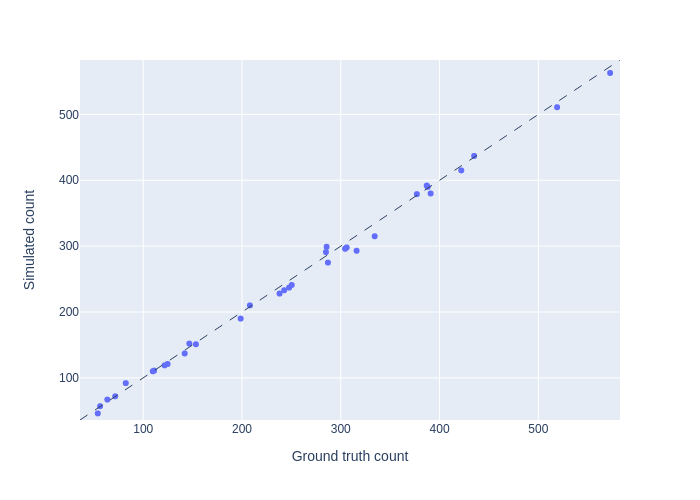}
         \caption{10:00 to 11:00}
         \label{fig:y equals x}
    \end{subfigure}
     \hfill
     \begin{subfigure}[b]{0.23\textwidth}
         \centering
         \includegraphics[width=\textwidth]{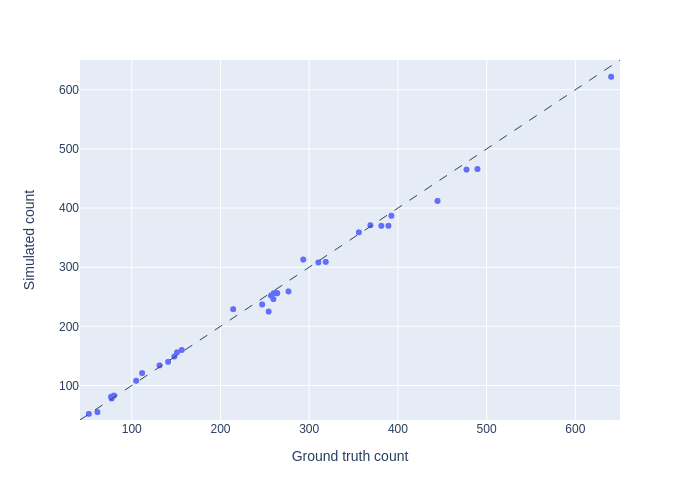}
         \caption{11:00 to 12:00}
         \label{fig:three sin x}
     \end{subfigure}
     \\
     \begin{subfigure}[b]{0.23\textwidth}
         \centering
         \includegraphics[width=\textwidth]{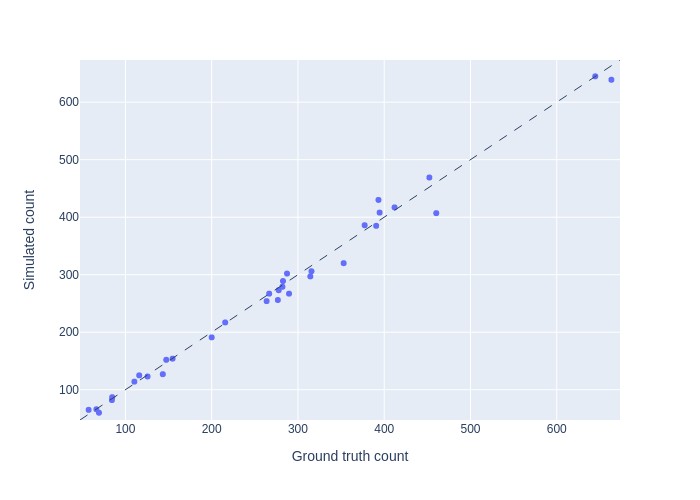}
         \caption{12:00 to 13:00}
         \label{fig:y equals x}
     \end{subfigure}
     \hfill
    \begin{subfigure}[b]{0.23\textwidth}
         \centering
         \includegraphics[width=\textwidth]{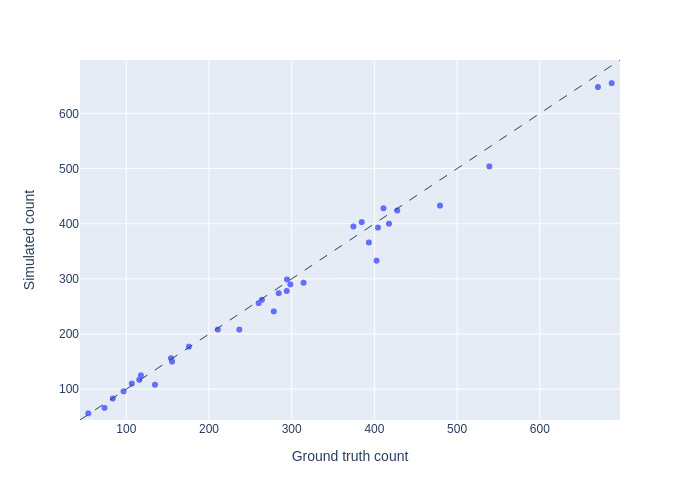}
         \caption{13:00 to 14:00}
         \label{fig:y equals x}
    \end{subfigure}
     \hfill
    \begin{subfigure}[b]{0.23\textwidth}
         \centering
         \includegraphics[width=\textwidth]{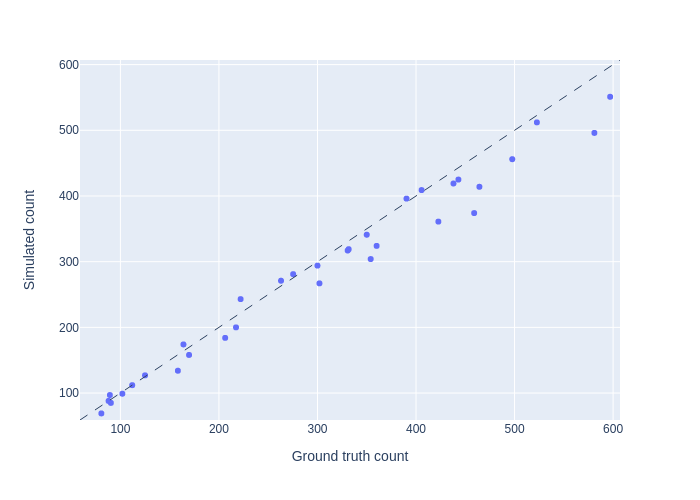}
         \caption{14:00 to 15:00}
         \label{fig:y equals x}
    \end{subfigure}
     \hfill
     \begin{subfigure}[b]{0.23\textwidth}
         \centering
         \includegraphics[width=\textwidth]{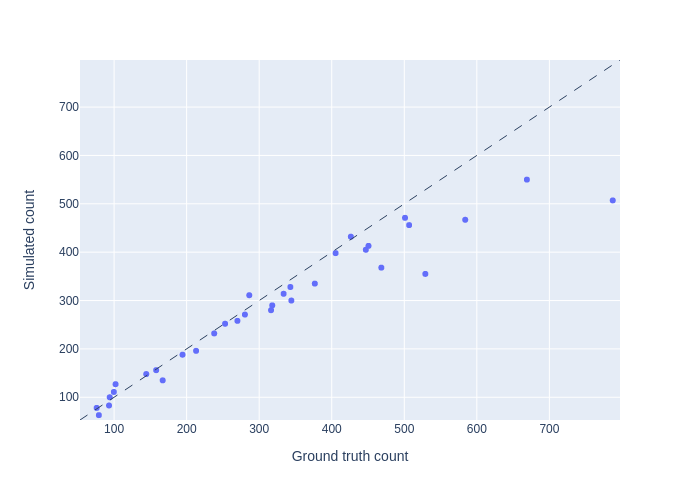}
         \caption{15:00 to 16:00}
         \label{fig:three sin x}
     \end{subfigure}
     \\
     \begin{subfigure}[b]{0.23\textwidth}
         \centering
         \includegraphics[width=\textwidth]{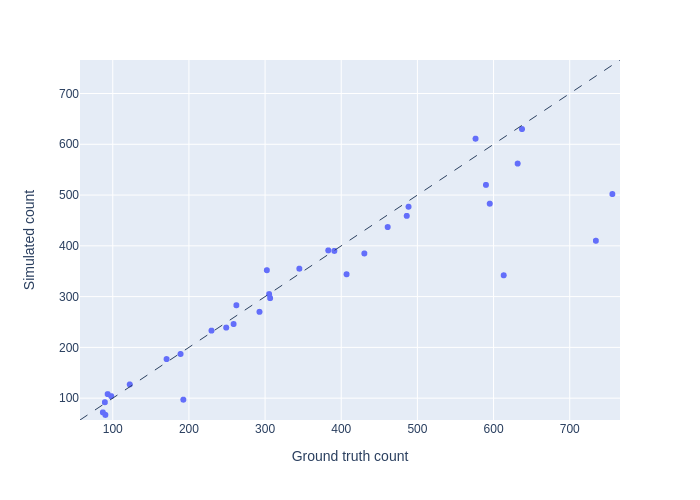}
         \caption{16:00 to 17:00}
         \label{fig:y equals x}
     \end{subfigure}
     \hfill
    \begin{subfigure}[b]{0.23\textwidth}
         \centering
         \includegraphics[width=\textwidth]{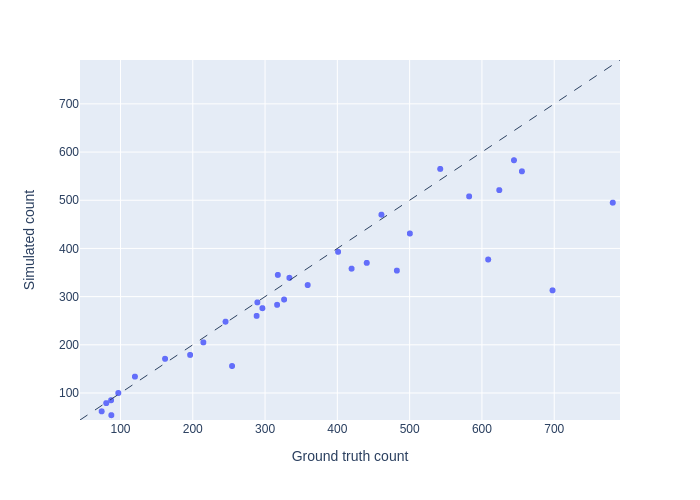}
         \caption{17:00 to 18:00}
         \label{fig:y equals x}
    \end{subfigure}
     \hfill
    \begin{subfigure}[b]{0.23\textwidth}
         \centering
         \includegraphics[width=\textwidth]{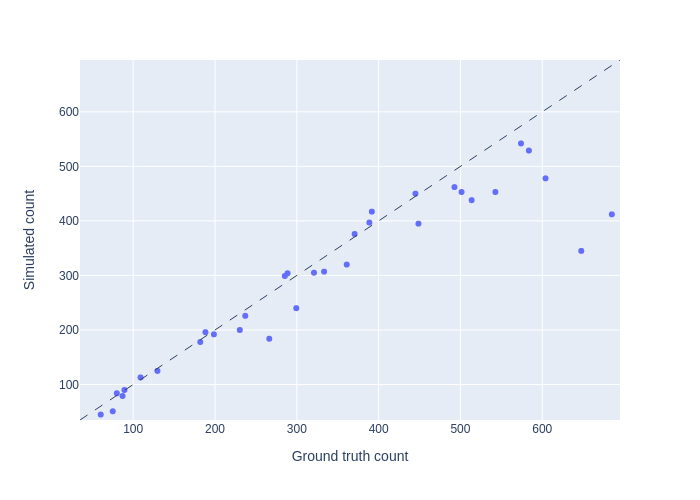}
         \caption{18:00 to 19:00}
         \label{fig:y equals x}
    \end{subfigure}
     \hfill
     \begin{subfigure}[b]{0.23\textwidth}
         \centering
         \includegraphics[width=\textwidth]{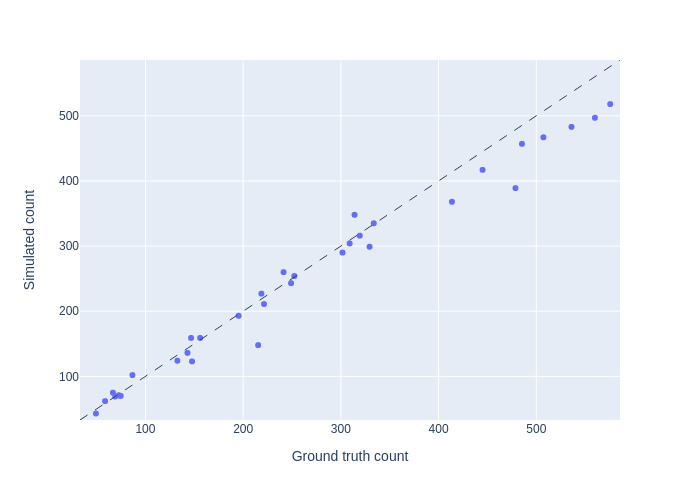}
         \caption{19:00 to 20:00}
         \label{fig:three sin x}
     \end{subfigure}
     \\
     \begin{subfigure}[b]{0.23\textwidth}
         \centering
         \includegraphics[width=\textwidth]{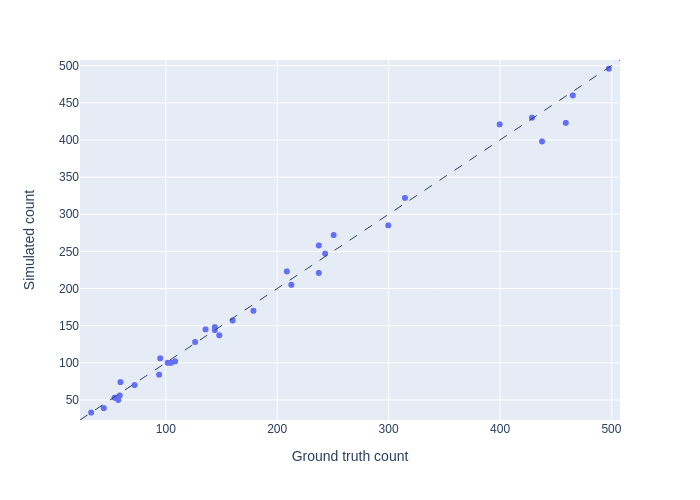}
         \caption{20:00 to 21:00}
         \label{fig:y equals x}
     \end{subfigure}
     \hfill
    \begin{subfigure}[b]{0.23\textwidth}
         \centering
         \includegraphics[width=\textwidth]{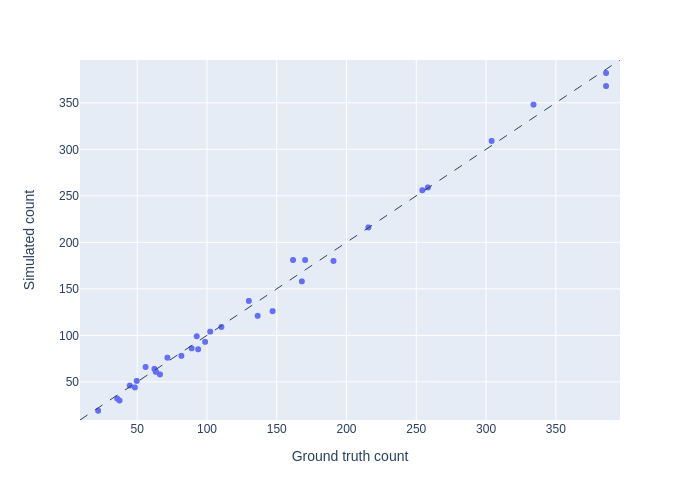}
         \caption{21:00 to 22:00}
         \label{fig:y equals x}
    \end{subfigure}
     \hfill
    \begin{subfigure}[b]{0.23\textwidth}
         \centering
         \includegraphics[width=\textwidth]{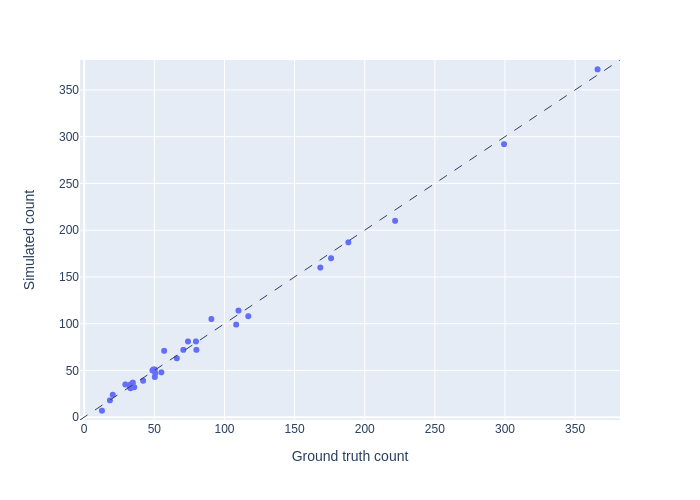}
         \caption{22:00 to 23:00}
         \label{fig:y equals x}
    \end{subfigure}
     \hfill
     \begin{subfigure}[b]{0.23\textwidth}
         \centering
         \includegraphics[width=\textwidth]{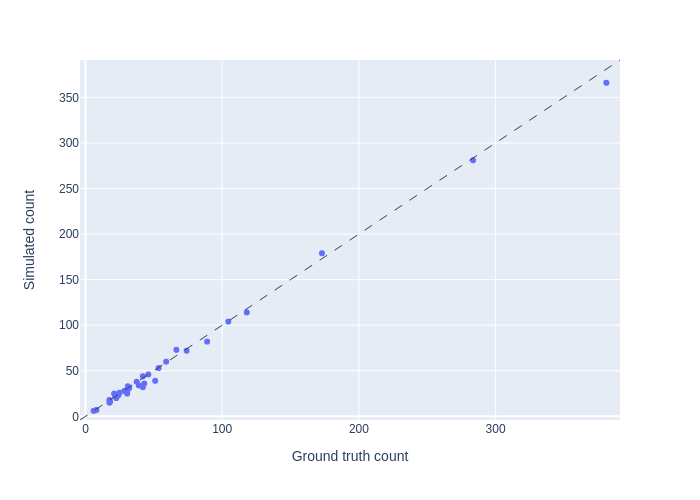}
         \caption{23:00 to 00:00}
         \label{fig:three sin x}
     \end{subfigure}
     \\
        \caption{Ground truth vs. simulated link count fits for different one-hour time frames of a sample day.  }
        \label{fig:24-hour}
\end{figure}

Figure \ref{fig:24-hour} compares the simulated link traffic counts and the ground truth counts within 24 hours of the experiment. The details of error calculations for each one-hour time frame are provided in table \ref{tab:24-hours}. Also, in addition to the number of estimated trips and average speed, table \ref{tab:24-hours} indicates the running time for each time frame and the number of calibration rounds before meeting the exit condition. It can be observed that although the limited computational budget of 10 calibration rounds is met in some time frames, the error rates are still reasonable. In peak hours, the high number of trips and the process of loading the network state have resulted in relatively higher error rates. While, there is always a trade-off between the computational budget and achieving more accurate results, we decided to limit the budget to stay relevant to the real-time application. The running time of the method for each one-hour time frame is, on average less than half an hour, with a maximum of 75 minutes.

\begin{table}[pos=h]
\centering
\caption{Hourly error estimations for the experiment in the city of Tartu with the running time and number of iterations of the calibration process. }
\begin{tabular}{c|ccccccc}
Time frame & \# trips & Average speed & $\epsilon$ \% & RMSE   & NRMSE \%   & Running time (sec)    &\# Iterations  \\ [1ex]
 \hline \\
0          & 1075            & 46.69        & 8.26      & 2.48   & 14.04    &719    &4     \\ [1ex]
1          & 693             & 46.58        & 7.19      & 1.75   & 11.72    &390    &1    \\ [1ex]
2          & 581             & 46.80        & 9.95      & 1.76   & 16.01    &2313   &7  \\ [1ex]
3          & 662             & 46.76        & 12.01     & 2.05   & 17.42    &3271   &10 \\ [1ex]
4          & 1063            & 47.55        & 9.46      & 2.96   & 12.61    &939    &3    \\ [1ex]
5          & 2272            & 46.83        & 7.42      & 4.55   & 9.02     &182    &1     \\ [1ex]
6          & 7239            & 44.82         & 3.45      & 6.06   & 4.08    &274    &1    \\ [1ex]
7          & 12934           & 37.54        & 9.69      & 30.79  & 11.46    &673    &2  \\ [1ex]
8          & 14728           & 33.04        & 12.74     & 48.72  & 14.85    &4116   &10  \\ [1ex]
9          & 12092           & 32.86        & 13.98     & 49.68  & 16.11    &3855   &10  \\ [1ex]
10         & 9234            & 39.42         & 3.05      & 8.67   & 3.48    &335    &1  \\ [1ex]
11         & 11529           & 42.44        & 4.29      & 12.55  & 4.90     &314    &1   \\ [1ex]
12         & 12070           & 41.25        & 5.24      & 16.33  & 6.01     &334    &1  \\ [1ex]
13         & 12996           & 38.23        & 6.62      & 22.03  & 7.58     &390    &1  \\ [1ex]
14         & 12529           & 37.33        & 9.45      & 31.80   & 10.59    &385   &1\\ [1ex]
15         & 13102           & 35.82         & 19.23     & 70.35  & 22.00    &3891  &10       \\ [1ex]
16         & 15002           & 32.65        & 22.79     & 93.50   & 26.01   &4304   &10  \\ [1ex]
17         & 12924           & 31.32         & 25.25     & 104.47 & 28.76   &4528   &10 \\ [1ex]
18         & 12413           & 34.02         & 21.59     & 82.05  & 24.59    &4105  &10\\ [1ex]
19         & 10114           & 35.46         & 10.11     & 31.17  & 11.82    &3612  &10   \\ [1ex]
20         & 8264            & 38.73        & 5.53      & 13.44  & 6.71     &300    &1   \\ [1ex]
21         & 6228            & 42.58        & 5.08      & 8.64   & 6.29    &256     &1 \\ [1ex]
22         & 4076            & 43.63        & 5.22      & 6.23   & 7.08     &298    &1\\[1ex]
23         & 3297            & 44.49        & 4.92      & 4.89   & 7.82    &385  &1
\end{tabular}
\label{tab:24-hours}
\end{table}

\begin{figure}[pos=h] 
\centering
\includegraphics[width=.8\textwidth]{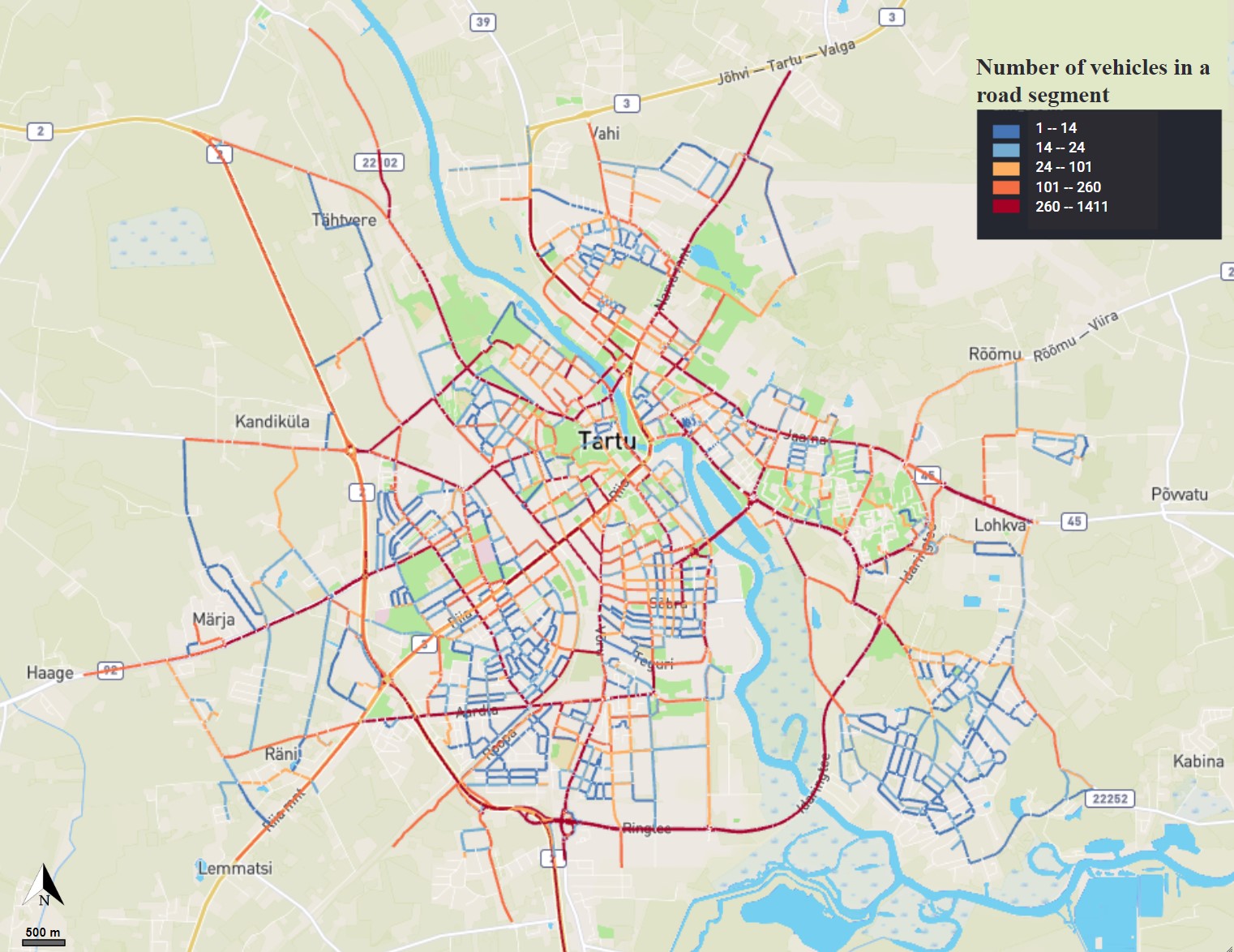}
\caption{The traffic heat-map of simulated trips at time frame 16:00 to 17:00. }
\label{fig:heatmap}
\end{figure}

A traffic heat map of the simulated trips for the time frame 16:00 to 17:00 is presented in figure \ref{fig:heatmap}, indicating the coverage of routes exploited by the generated trips of the method. 

\newpage
\section{Conclusions} \label{sec:conclusion}
This paper presents a computationally efficient algorithm for high-dimensional real-time dynamic demand estimation in congested networks. On a high level, the method follows the simulation-based bi-level optimization framework. Unlike typical approaches, the algorithm estimates the demand sequentially by splitting the computations into discrete concise time frames. To ensure the temporal dependency of the traffic status in consecutive time frames, the microsimulation of the network at the end of each time frame is loaded to the next time frame as the initial network state. Hence, no trip is left incomplete, and the micro-level transitions between time frames are as smooth as possible. The sequential approach allows the method to receive the field measurements input as a data stream, which makes the solution compatible with real-time and online applications. Most previous works in the literature rely on a prior OD matrix (seed) to estimate accurately. In contrast, this work reduces the input required to the easy-to-access normalized OD distribution that can be obtained from various available sources, such as population distribution. 

The upper-level optimization problem is defined as a bounded variable quadratic programming problem. Since the coefficients' matrix is relatively sparse, it is possible to take advantage of efficient off-the-shelf algorithms for optimizing this problem. The lower-level DTA problem is modeled using probabilistic and stochastic route choice methods. The probabilistic parameters of route assignment and trip generations on the route level are estimated through a feedback loop between supply and demand models by iteratively solving the upper-level optimization problem and performing parallel samplings and simulations for verification. This approach avoids costly black-box simulations to solve the DTA problem. A disaggregated route-level demand estimation is obtained at the end of the calculations in each time frame. Moreover, we can avoid rounding errors in converting the real-valued outcomes to integers with the probabilistic approach. In addition, the convergence is accelerated by formulating the problem as a fixed-point problem for estimating the link travel times of the network for each time frame. 

The method is validated in synthetic and real-world settings, with detailed error analysis. The case study of the city of Tartu presents the application of the technique to a high-dimensional large-scale network with a stream of real traffic counts collected by IoT sensors. The hourly simulation outcomes show a stable performance with high accuracy under a tight computational budget for a whole day. Indeed, an additional computational budget could improve the method's precision, especially during peak hours. 

In future works, one can include dynamic NOD for different hours of the day instead of a static NOD input. The route choice can be modeled more realistically than the classic logit model. Moreover, in practice, by embedding more spatial information in selecting the points of interest as the OD pairs, the method will be able to model the real dynamics of a city more accurately.  
The micro-level output of the presented work can provide valuable input for prediction models and machine learning approaches.


\section*{Acknowledgment}
This work was supported by the European Social Fund via IT Academy programme, and the Estonian Centre of Excellence in IT (EXCITE).

\newpage
 \appendix
 \section{Notations}
 The following is the list of notations used in the paper. 

\begin{table*}[pos=hbt!]
\begin{tabular}{l l}
$\mathcal{N}$ & Set of nodes in network; \\
$\mathcal{L}$ & Set of inks in the network; \\
$\mathcal{Q}$ & Set of links with traffic counter;\\
$q$ & The number of traffic counters in the network;\\
$W$ & The set of OD pairs;\\
$|W|$ & The number of OD pairs;\\
$t\in \{0, 1, 2, \cdots\}$ & Time frame index;\\
$i$ & Route index;\\
$m$ & OD pair index;\\
$k$ & Traffic counter (sensor) index;\\
$\Delta$ & Time frame length;\\
$\delta$ & Departure time;\\
$\mathbf{\Tilde c}=(\Tilde c_k)$ & Vector of field traffic counts; \\
$\mathbf{c}^* = (c^*_k)$ & Vector of estimated traffic counts; \\
${\bf N} = (\eta_{m})$ & Normalized OD trip distribution;\\
${\bf \Tilde X}$ & Prior OD matrix (seed); \\
$\sigma$ & Scaling factor in ${\bf \Tilde X} = \sigma \times {\bf N} $;\\
${\bf X^*}=(x^*_m)$ & Expected OD matrix; \\
$R_{m}$ & Set of routes between $m$th OD pair;\\
$\theta_{i}$ & Travel time of route $i$;\\
$\theta_{ik}$ & Travel time from origin to sensor $k$, on route $i$; \\ 
$E_i$ & Expected number of trips on route $i$;\\
${\bf A} = (\alpha_{mk}) $ & Assignment matrix;  $\alpha_{mk}$ is the probability of a trip between $m$th OD pair crosses sensor $k$;\\
${\bf \tau}$ & Vector of link travel times for a fixed time frame. \\

\end{tabular}
\end{table*}



\bibliographystyle{cas-model2-names}
\bibliography{cas-refs}





\end{document}